\RequirePackage{fix-cm}
\documentclass[smallextended]{svjour3}       
\smartqed  
\usepackage{graphicx}
\usepackage{amssymb,amsmath,graphicx,color,pagesize}
\usepackage{amsfonts}
\usepackage[caption=false]{subfig}
\usepackage{epsfig,latexsym,graphicx}
\begin{document}

\title{The Minimum $S$-Divergence Estimator under Continuous Models: The Basu-Lindsay Approach
}

\titlerunning{The MSDE under Continuous Models: The Basu-Lindsay Approach}        

\author{Abhik Ghosh         \and
	Ayanendranath Basu 
}


\institute{Abhik Ghosh  \at
	Indian Statistical Institute, Kolkata, India \\
	\email{abhianik@gmail.com}           
	\and
	Ayanendranath Basu \at
	Indian Statistical Institute, Kolkata, India \\
	\email{ayanendranath.basu@gmail.com}  
}
\date{Received: date / Accepted: date}

\maketitle

\begin{abstract}
	Robust inference based on the minimization of  statistical divergences has proved to be a useful 
	alternative to the classical maximum likelihood based techniques. Recently Ghosh et al.~(2013) 
	proposed a general class of divergence measures for robust statistical inference, named the
	$S$-Divergence Family. Ghosh (2014) discussed its asymptotic properties for 
	the discrete model of densities. 
	In the present paper, we develop the asymptotic properties of the proposed minimum $S$-divergence 
	estimators under continuous models. Here we use the Basu--Lindsay approach (1994) of smoothing 
	the model densities that, unlike previous approaches, avoids much of the complications of the 
	kernel bandwidth selection. Illustrations are presented to support the performance of the 
	resulting estimators both in terms of efficiency and robustness through extensive simulation studies 
	and real data examples. 
\end{abstract}

\section{Introduction}

In case of parametric statistical inference, the usefulness and properties of various
density based minimum divergence estimators have been extensively studied in the recent literature; 
see Basu et al.~(2011). 
The key in density-based minimum divergence estimation is the quantification of the discrepancy between
the parametric model and the sample data through a suitable density based divergence. 
Then one obtains the estimate of the unknown parameter of interest
by minimizing this divergence as a function of the unknown parameter. 
The class of density-based divergences that may be 
useful in this approach includes the Pearson's chi-square (Pearson, 1900), 
the Kullback-Leibler divergence (Kullback and Leibler, 1951), 
the $\phi$-divergence or disparity family (Csisz\'{a}r, 1963) including the Hellinger distance, 
the Bregman divergence (Bregman, 1967), the Burbea-Rao divergence (Burbea and Rao, 1982), 
Cressie-Read family of Power divergence (Cressie and Read, 1984), 
the density power divergence (Basu et al., 1998) etc.
The power divergence family is a subclass of the class of $\phi$-divergences 
while the density power divergence family is a subclass of Bregman divergences.  
There are many examples of the application of the minimum divergence estimation method in recent literature; 
see, for example, Lindsay (1994), Menendez  et al.~(1995, 1998), 
Jones et al.~(2001), Mihoko and Eguchi (2002), Landaburu et al.~(2005), Pardo et al.~(2006), 
Martin and Pardo (2008), Lee and Lee (2009), Menendez, Pardo and Pardo (2009), Toma and Broniatowski (2010),
Ghosh and Basu (2013a,b, 2014a,b) among many others.

The main advantage of minimum distance methods in parametric estimation 
is the robustness property which some of them inherently possess. 
In the past statisticians believed that the goals of 
robustness and efficiency were conflicting and could not be achieved simultaneously.  
However, several density-based minimum distance estimators have been demonstrated 
to have strong robustness properties along with full asymptotic efficiency 
[e.g., Beran (1977), Tamura and Boos (1986), Simpson (1987, 1989) and Lindsay (1994)]. 
Several of the minimum  distance estimators based on the Cressie-Read family of power divergences are 
examples of such estimators. However, standard estimation techniques under continuous models using 
the power divergence family require the use of nonparametric smoothing for the 
construction of the data generating density and hence inherit all the complications 
of the kernel density estimation process and bandwidth selection. 
Later, Basu et al.~(1998) developed a rich class of density-based divergence measures 
named the density power divergences that produce robust estimates without using nonparametric smoothing
but with a small loss in efficiency. Both the density power divergence and the power divergence families have 
similarities in their outlier downweighting philosophy. Although these two families have different
forms (except for the Kullback-Leibler divergence which is the only divergence common to both families), 
Patra et al.~(2013) provided an useful connection between these
two families which allows us to develop either family from the other through specific motivations.

Combining the forms of the power divergence and the density power divergences, 
Ghosh et al.~(2013) developed a two parameter family of density-based divergences, called the ``$S$-Divergence", 
that connects each member of the Cressie-Read family of power divergences smoothly to 
the $L_2$-divergence at the other end. Through various numerical examples, 
these authors illustrated that several of the minimum divergence estimators within the $S$-divergence family 
also have strong robustness properties and are often competitive with the classical estimators 
in terms of efficiency. 
Ghosh (2014) has provided the asymptotic distribution of the corresponding minimum $S$-divergence estimators 
under discrete model families.
In the discrete case there is a natural nonparametric density estimator of the true unknown density; 
this is simply the vector of the sample relative frequencies, the construction of which requires 
no additional artifacts like the kernel function or the  bandwidth. 
However, in the case of continuous models it is not so simple and one needs to use 
the kernel density estimate for the same purpose.

In this present article, we will develop the theoretical properties of the 
minimum $S$-divergence estimators under the set up for a continuous model. 
The first approach under the similar set-up with kernel density estimators was  presented by Beran (1977)
in the context of minimum Hellinger distance estimators. But, his approach 
contains all the complications of the kernel estimator along with the issue of bandwidth selection. 
For simplicity, we will  use the later approach of Basu and Lindsay (1994) that 
helps us to at least partially avoid the problem with  kernel bandwidth selection. 
In this approach we  replace ordinary model densities with the  smoothed model densities 
to obtain an estimator slightly different than the minimum $S$-divergence estimator; 
we will refer to this estimator as the minimum $S^*$-divergence estimator.
We will prove the  consistency and asymptotic normality of the minimum $S^*$-divergence estimator 
under appropriate assumptions. 
Interestingly, we will see that the asymptotic
distributions of the minimum $S^*$-divergence estimators are also independent of the parameter
$\lambda$ in the definition of $S$-divergence, 
as in the distribution of the minimum $S$-divergence estimators in the discrete models. 

The rest of the paper is organized as follows.
We start with a brief description of the $S$-divergence family and the corresponding 
minimum $S$-divergence estimators in Section \ref{SEC:MSDE}. 
Then we will describe the two different approaches of the minimum divergence estimator under continuous models, 
Beran's approach and Basu-Lindsay approach, in Section \ref{SEC:MS*DE}. 
The minimum $S^*$-divergence estimators will be introduced in that section and details 
for the particular case of $\lambda=0$ will be presented in Section \ref{SEC:MS*DE0}. 
Sections \ref{SEC:MS*DE_IF} and \ref{SEC:MS*DE_asymptotic} will describe the influence function analysis 
and the asymptotic properties of the general minimum $S^*$-divergence estimators. 
All the results will be supported by suitable simulation studies in Section \ref{SEC:4MS*DE_simulation} 
and some interesting real data examples in Section \ref{SEC:MS*DE_real data}. 
We will introduce the concept of the $\alpha$-transparent kernels in Section \ref{SEC:4Trans_kernel}. 
We will end the paper with a small discussion on the choice of tuning parameters in Section \ref{SEC:discussion}.

\section{The $S$-Divergence Family and its Estimating Equation}\label{SEC:MSDE}

The $S$-divergence family of divergences (Ghosh et al.~2013) is defined in terms of two parameters $\alpha \geq 0$
and $\lambda \in \mathbb{R}$ as
\begin{eqnarray}
S_{(\alpha,\lambda)}(g,f)
&=&  \frac{1}{A}~\int~f^{1+\alpha} - \frac{1+\alpha}{AB}~\int~~f^{B} g^{A} + \frac{1}{B}~\int~~g^{1+\alpha}, 
\label{EQ:S_div}
\end{eqnarray}
with $ A = 1+\lambda (1-\alpha)$ and $B = \alpha - \lambda (1-\alpha)$.
Note that, $ A+B=1+\alpha $. If $A=0$ then the corresponding $S$-divergence measures are defined 
as the continuous limit of (\ref{EQ:S_div}) as $A \rightarrow 0$ and are given by
\begin{eqnarray}
S_{(\alpha,\lambda : A = 0)}(g,f) &=& \lim_{A \rightarrow 0} ~ S_{(\alpha, \lambda)}(g,f) \nonumber\\
&=&  \int f^{1+\alpha} \log\left(\frac{f}{g}\right) - \int \frac{(f^{1+\alpha} - g^{1+\alpha})}{{1+\alpha}}.
\label{EQ:S_div_A0}
\end{eqnarray}
Similarly, for $B=0$ then the corresponding $S$-divergence measures are defined as
\begin{eqnarray}
S_{(\alpha,\lambda : B = 0)}(g,f) &=& \lim_{B \rightarrow 0} ~ S_{(\alpha, \lambda)}(g,f) \nonumber\\
&=&  \int g^{1+\alpha} \log\left(\frac{g}{f}\right) - \int \frac{(g^{1+\alpha} - f^{1+\alpha})}{{1+\alpha}}.
\label{EQ:S_div_B0}
\end{eqnarray}
Note that for $\alpha = 0$, this family reduces to the Cressie-Read family with parameter
$\lambda$ and for $\alpha =1$, it gives the $L_2$ divergence (independently of $\lambda$). 
Again for $\lambda = 0$, the expression in (\ref{EQ:S_div}) generates the Density Power divergence measure having parameter $\alpha$. 
Ghosh et.~al. (2013) showed that the $S$-divergence family defined in (\ref{EQ:S_div}),
(\ref{EQ:S_div_A0}) and (\ref{EQ:S_div_B0}) indeed represents a family of genuine statistical 
divergence measures in the sense that $S_{(\alpha, \lambda)}(g, f) \ge 0$ for all densities 
$g, f $ defined with respect to a common measure and for all $\lambda \in \mathbb{R}$, $\alpha \ge 0$, 
and it equals zero if and only if the argument densities are identically equal.


Now, let us assume that we have a sample from the true unknown density $g$ 
(having $G$ as the corresponding cumulative distribution function), and 
we want to model it by the parametric family $\{f_\theta : \theta\in\Theta\subseteq\mathbb{R}^p\}$ of densities.
The minimum $S$-divergence estimator of the unknown parameter $\theta$ is obtained by 
minimizing the divergence $S_{(\alpha, \lambda)}(\hat{g},f_{\theta})$ with respect to $\theta\in\Theta$,
where $\hat{g}$ is some nonparametric estimate of the true density $g$ based on the observed data. 
For the discrete model, the vector of relative frequencies may be taken as $\hat{g}$, 
but for continuous models we need to use kernel density estimation or a similar technique to obtain $\hat{g}$. 
The corresponding estimating equation for the minimum $S$-divergence estimator is given by
\begin{eqnarray}
\int f_{\theta}^{1+\alpha} u_{\theta} - \int f_{\theta}^{B} \hat{g}^{A} u_{\theta} &=& 0, \nonumber \\
\mbox{or,   } ~~ \int K(\delta(x))f_{\theta}^{1+\alpha}(x) u_{\theta}(x) &=& 0, \label{EQ:S-divergence_est_equation}
\end{eqnarray}
where $\delta(x) = \frac{\hat{g}(x)}{f_{\theta}(x)}-1$, $ K(\delta) = \frac{(\delta+1)^A - 1}{A} $,
$u_\theta(x) = \frac{\partial}{\partial\theta} \log f_\theta(x)$, the maximum likelihood score function.

Let $T_{\alpha,\lambda}$ denotes the minimum $S$-divergence functional at 
the distribution $G$ defined by the relation
$$
S_{(\alpha,\lambda)}(g, f_{T_{\alpha,\lambda}(G)}) = \min_{\theta \in \Theta} S_{(\alpha,\lambda)}(g, f_\theta),
$$
whenever the minimum exists. 
One can easily derive the influence function of this functional by taking a derivative of 
the estimating equation (\ref{EQ:S-divergence_est_equation}) under a mixture contaminated density 
with respect to a mixture contamination. 
Then, the influence function of this minimum $S$-divergence functional is given by
\begin{eqnarray}
IF(y;T_{\alpha,\lambda},G) = J^{-1} \left[ A u_{\theta}(y) f_{\theta}^B(y) g^{A-1}(y) - \xi \right],
\end{eqnarray}
where 
\begin{eqnarray}
\xi = \xi(\theta) &=& A \int u_{\theta}f_{\theta}^B g^{A},\\
J = J(\theta)  &=& A \int u_{\theta}u_\theta^T f_{\theta}^{1+\alpha}  +  
\int (i_{\theta} - B u_{\theta}u_\theta^T)(g^A - f_{\theta}^A) f_{\theta}^B,
\end{eqnarray}
and $i_{\theta}(x) = -\nabla[u_{\theta}(x)]$ with $\nabla$ representing the gradient with respect to $\theta$. 
However, if the true density belongs to the model family, i.e., $g=f_{\theta}$, 
the influence function simplifies to 
\begin{eqnarray}
IF(y;T_{\alpha,\lambda},F_\theta) =  \left(\int u_{\theta}u_\theta^T f_{\theta}^{1+\alpha}\right)^{-1} 
\left[u_{\theta}(y) f_{\theta}^{\alpha}(y)  - \int u_{\theta}f_{\theta}^{1+\alpha}\right].
\label{EQ:S_div_IF_model}
\end{eqnarray}
Interestingly the influence function at the model depends only on the parameter $\alpha$ 
and not on $\lambda$.

\section{The Minimum $S$-Divergence Estimator under Continuous Models: Different Approaches}\label{SEC:MS*DE}

Let us now consider the estimation under the continuous models by minimizing the $S$-divergence. 
Let $\mathcal{G}$ denotes the class of all probability distributions having densities 
with respect to the Lebesgue measure. We will assume that the true, data generating distribution 
$G$ and the model family $\mathcal{F} = \{ F_\theta : \theta \in \Theta \subseteq \mathbb{R}^p \}$ belong to $\mathcal{G}$;
also $G$ and $F_\theta$ have densities $g$ and $f_\theta$ with respect to the Lebesgue measure.
\index{MSDE!under Continuous Models}

Let $X_1,\ldots,X_n$ be a random sample of size $n$ from the true distribution $G$ which is 
modeled by $\mathcal{F}$, and we want to estimate the unknown model parameter $\theta$. 
As in the case of discrete models, the minimum $S$-divergence estimator (MSDE) of 
the unknown parameter $\theta$ is to be obtained  by choosing the model density $f_\theta$ 
which gives the closest fit to the data with respect to the $S$-divergence  measure. 
However, unlike the discrete case, this gives rise to an immediate challenge; 
the data are discrete, but the model is continuous, so that there is an obvious incompatibility of measures 
in constructing a distance between the two. We cannot simply use relative frequencies to represent a 
nonparametric density estimate of the true data generating distribution in this case.


In this context, Beran (1977) suggested the construction of a continuous density estimate using 
some appropriate nonparametric density estimation method such as kernel density estimation.
Let us assume that
\begin{equation}
g^*_n(x) = \frac{1}{n}  \sum_{i=1}^n  W (x, X_i, h_n) = \int 
W(x, y, h_n) dG_n (y)    
\label{EQ:g-star-n}
\end{equation}
represents a nonparametric kernel density estimator 
where $W(x, y, h_n)$ is a smooth kernel function with bandwidth $h_n$ and
$G_n$ is the empirical distribution function as obtained from the data. 
Usually the kernel is chosen as a symmetric location-scale density with scale $h_n$, i.e.,  
\begin{equation}
W(x, X_i, h_n) = \frac{1}{h_n} w\left(\frac{x - X_i}{h_n}\right),
\end{equation}
where $w(\cdot)$ is a symmetric nonnegative density function.
According to Beran's approach, we can now estimate $\theta$ by minimizing the $S$-divergence  measure 
between $g_n^*$ and $f_\theta$. Under suitable differentiability assumptions, the estimating 
equation is then given by (\ref{EQ:S-divergence_est_equation}) with $\hat{g}$ replaced by 
the kernel estimator $g_n^*$. 
Although  the rest of the estimation and the asymptotic properties of the estimator 
are similar to the discrete case, the inclusion of the kernel density estimation process leads to 
substantial difficulties in practice. The theoretical derivation of the asymptotic normality of the minimum 
$S$-divergence  estimators and the description of their other asymptotic properties are far more complex in 
this case. In particular, the choice of the sequence of kernels (or, more precisely, the sequence 
of smoothing parameters)  becomes critical and also various complicated conditions have to be imposed 
for the consistency of the kernel density estimator (see Wied and Weißbach, 2012, for a survey);
these complications are unavoidable for Beran's approach.

The approach taken by Basu and Lindsay (1994) differs from the Beran's approach in that 
it proposes that the model be convoluted with the same kernel as well. 
Suppose $f_\theta^*$ represents the kernel integrated ``smoothed" version of the model defined as
\index{Kernel Smoothed Model}
\begin{equation}
f^*_\theta(x) = \int W(x, y, h) dF_\theta(y).
\label{EQ:3smoothed_model}
\end{equation}
In the Basu-Lindsay approach of model smoothing, 
we consider the $S$-divergence  measure between $g_n^*$ and $f_\theta^*$, 
and minimize it over $\theta$ to obtain the corresponding minimum $S$-divergence estimator. 
The procedure may be intuitively justified as follows. 
The intent here is to minimize a measure of discrepancy between the data and the model.
To make the data continuous, one has to import an artificial kernel. However, one needs 
to ensure $-$ through the imposition of suitable conditions on the kernel function and the smoothing 
parameter $-$ that the additional smoothing effect due to the kernel vanishes asymptotically. 
In the Basu-Lindsay approach, we convolute the model 
with the same kernel used on the data. In a sense, this compensates for the bias due to the imposition
of the kernel on the data by imposing the same bias on the model. It is therefore expected that the
kernel will play a less important role in the estimation procedure than it plays in Beran's approach,
particularly in small samples. As we will see later in this section, one gets consistent estimators of 
the parameter $\theta$ even when the smoothing parameter is held fixed as the sample size 
increases to infinity.

In this work, we will consider only the simpler Basu-Lindsay approach for estimation 
under the continuous model based on $S$-divergence. We need consider the corresponding smoothed 
version of the likelihood score function 
$$\widetilde{u}_{\theta}(x) = \nabla \log f_\theta^*(x);$$ 
we will denote its $j^{\rm th}$ entry by $\widetilde{u}_{j\theta}(x)$.
Further, for the second derivative matrix $\nabla\widetilde{u}_\theta(x)=\nabla^2 \log f_\theta^*(x)$, 
we denote its $(j,k)^{\rm th}$ element by $\widetilde{u}_{jk\theta}(x)$.
Now, let us consider  $S_{(\alpha, \lambda)}(g_n^*, f_\theta^* )$ and minimize this measure of 
discrepancy between the model and data with respect to the parameter $\theta$.  
A routine differentiation shows that the corresponding estimating equation is given by 
\begin{eqnarray}
\int (f_\theta^*)^B(g_n^*)^A \widetilde{u_\theta} - \int (f_\theta^*)^{1+\alpha} \widetilde{u_\theta} &=& 0, \nonumber\\
\mbox{or,} ~~~~~~~ \int K(\delta_n^*(x)) (f_\theta^*)^{1+\alpha} \widetilde{u_\theta} &=& 0,
\label{EQ:S-div_est_eq_kernel} 
\end{eqnarray}
with $\delta_n^*(x) = \frac{g_n^*(x)}{f_\theta^*(x)} -1$.
We will denote the estimator obtained by this approach as minimum $S^*$-divergence estimator
($\textrm{MSDE}^*$) which is, in general, not the same as the estimator obtained by minimizing 
the $S$-divergence measure  $S_{(\alpha, \lambda)}(g_n^*, f_\theta )$.
\index{MSDE*!Minimum S*-divergence estimator} 
The main reason behind this inequality is the
substitution of the model family by its kernel smoothed version.

We have commented earlier that our method, based on the Basu-Lindsay approach, 
leads to consistent estimators for a fixed bandwidth $h$, which makes the bandwidth selection less critical.
While we will formally prove this in Section \ref{SEC:MS*DE_asymptotic}, 
here  we point out some practical advantages of this method.
In the ordinary MSDE (using Beran's approach), the data alone are smoothed with the kernel,
so that the scale of the data distribution gets inflated with increasing bandwidth without a 
corresponding compensation in the model distribution. 
As a result, the MSDE of the scale parameter (in the normal model, for example)
also gets inflated with increasing bandwidth. 
This phenomenon is either absent or significantly less prominent  for the MSDE$^*$ 
(in Basu -Lindsay approach), which underlines the fact that the scale parameters 
have little dependence on the bandwidth in the Basu-Lindsay approach.

We explore this issue further with a numerical example, 
for which we choose the pseudo-random sample $X_1, X_2, \ldots, X_{40}$ from $N(0,1)$ 
selected and analyzed by Beran (1977; Section 6, p.~460).
For these data, we determine the MSDE and MSDE$^*$s of $\mu$ and $\sigma$ 
under the normal $N(\mu, \sigma^2)$ model for several combinations of $\alpha$ and $\lambda$ 
at different values of the smoothing parameter. 
The iterative process uses the initial values $\mu_0 = \mbox{median}_i \{X_i\}$ for the mean and  
\begin{equation}
{\sigma_0} = \frac{{\rm median}_i {|X_i - \mu_0| }}{0.6745}
\label{EQ:sigma_0_est}
\end{equation}
for the scale. Several values of the bandwidth $h$ are chosen as $h=h_0\sigma_0$,
where $h_0$ varies from $0.4$ to $1$ at intervals of $0.1$. 
The scale estimates for each choice of $\alpha, \lambda, h$ and 
the estimation method are reported in Table \ref{TAB:sigma_comp}.

 
\begin{table}[!th]
	\centering 
	\caption{Estimates of the Scale parameter $\sigma$ for Beran's dataset }
	\begin{tabular}{r|ll|lllllll} \hline\noalign{\smallskip}
	&	$\lambda$	&	$\alpha$	&		&		&		&	$h_0$	&		&		&		\\
	&		&		&	0.4	&	0.5	&	0.6	&	0.7	&	0.8	&	0.9	&	1	\\
	\noalign{\smallskip}\hline\noalign{\smallskip}
MSDE	&	0	&	0.5	&	1.0298	&	1.0705	&	1.1166	&	1.1674	&	1.2223	&	1.2807	&	1.3423	\\
	&	0	&	0.3	&	1.0493	&	1.0860	&	1.1287	&	1.1767	&	1.2294	&	1.2861	&	1.3464	\\
	&	$-0.5$	&	0	&	1.0480	&	1.0877	&	1.1315	&	1.1798	&	1.2321	&	1.2885	&	1.3483	\\
	&	$-0.5$	&	0.3	&	1.0405	&	1.0803	&	1.1249	&	1.1741	&	1.2275	&	1.2848	&	1.3455	\\
	&	$-$0.5	&	0.5	&	1.0233	&	1.0663	&	1.1140	&	1.1657	&	1.2210	&	1.2798	&	1.3417	\\
	&	$-$1	&	0.5	&	1.0165	&	1.0621	&	1.1112	&	1.1638	&	1.2198	&	1.2790	&	1.3410	\\
	&	$-$1	&	0.3	&	1.0301	&	1.0739	&	1.1209	&	1.1714	&	1.2257	&	1.2834	&	1.3444	\\
	&	1	&	0.5	&	1.0415	&	1.0781	&	1.1218	&	1.1709	&	1.2247	&	1.2824	&	1.3435	\\
	&	1	&	0	&	1.0810	&	1.1099	&	1.1469	&	1.1906	&	1.2401	&	1.2944	&	1.3529	\\
	&	-	&	1	&	0.9809	&	1.0322	&	1.0874	&	1.1454	&	1.2059	&	1.2685	&	1.3332	\\
	\noalign{\smallskip}\hline\noalign{\smallskip}
MSDE$^*$	&	0	&	0.5	&	0.9637	&	0.9694	&	0.9747	&	0.9792	&	0.9830	&	0.9860	&	0.9885	\\
	&	0	&	0.3	&	0.9844	&	0.9865	&	0.9885	&	0.9902	&	0.9918	&	0.9930	&	0.9941	\\
	&	$-$0.5	&	0	&	0.9831	&	0.9884	&	0.9917	&	0.9938	&	0.9952	&	0.9961	&	0.9968	\\
	&	$-$0.5	&	0.3	&	0.9750	&	0.9803	&	0.9842	&	0.9872	&	0.9896	&	0.9914	&	0.9928	\\
	&	$-$0.5	&	0.5	&	0.9567	&	0.9648	&	0.9716	&	0.9771	&	0.9814	&	0.9849	&	0.9877	\\
	&	$-$1	&	0.5	&	0.9494	&	0.9602	&	0.9685	&	0.9749	&	0.9799	&	0.9838	&	0.9868	\\
	&	$-$1	&	0.3	&	0.9639	&	0.9733	&	0.9795	&	0.9840	&	0.9872	&	0.9897	&	0.9915	\\
	&	1	&	0.5	&	0.9760	&	0.9779	&	0.9806	&	0.9833	&	0.9860	&	0.9882	&	0.9902	\\
	&	1	&	0	&	1.0180	&	1.0128	&	1.0093	&	1.0068	&	1.0050	&	1.0038	&	1.0028	\\
	&	-	&	1	&	0.9112	&	0.9270	&	0.9411	&	0.9529	&	0.9625	&	0.9702	&	0.9761	\\
		\noalign{\smallskip}\hline
	\end{tabular}
	\label{TAB:sigma_comp}
\end{table}

The results show that the estimates of the scale parameter $\sigma$ 
are remarkably stable over variation in $h_0$ when one uses the MSDE$^*$, but this is not the case for the MSDE.
Actual calculation shows that, depending on the $(\alpha, \lambda)$ combination, 
the variability in the MSDE of scale over the range of $h_0$ considered in this example
is between 6 to 30 times the corresponding variability in the MSDE$^*$.  
On the other hand the location estimates for the two methods (not presented here for the sake of brevity)
are extremely close and show very little bandwidth effect.
This is intuitively expected as the effect of smoothing is primarily on the scale parameter 
and not on the location parameter.
This example shows that the Basu-Lindsay approach not only makes the choice of a proper bandwidth 
theoretically unimportant, it also substantially reduces the impact of the bandwidth 
on the estimates of the scale parameter in finite samples.

\section{The Minimum $S^*$-Divergence Estimator: Special Case ($\lambda = 0$)}\label{SEC:MS*DE0}

The density power divergence measure, originally proposed by Basu et al.~(1998), can also be obtained 
from the $S$-divergence divergence family by putting $\lambda = 0$ and has the form
\begin{eqnarray}
d_\alpha(g,f_\theta) &=& S_{(\alpha, 0)}(g,f_\theta) \nonumber\\ 
&=&  \int ~ f_\theta^{1+\alpha}  -   \frac{1+\alpha}{\alpha} ~ \int ~~ f_\theta^{\alpha} g  
+ \frac{1}{\alpha} ~ \int ~~ g^{1+\alpha}, ~~ \alpha>0. \nonumber
\end{eqnarray}
This particular family of divergences has the special property that we can find the 
minimum DPD estimator without using kernel density estimators even under continuous models. 
In fact, we can write the minimum DPD estimating equation in terms of the empirical distribution 
function $G_n$ based on the observed data as 
\begin{eqnarray}
\int f_\theta^\alpha {u_\theta} dG_n - \int f_\theta^{1+\alpha} {u_\theta} &=& 0, \nonumber\\
\mbox{or,} ~~~ 
\frac{1}{n}  \sum_{i=1}^n f_\theta^\alpha(X_i) {u_\theta}(X_i) - 
E_\theta [ f_\theta^{\alpha}(X) {u_\theta}(X)] &=& 0. \label{EQ:est-eqn-MDPD}
\end{eqnarray}
Note that the above estimating equation is an unbiased estimating equation at the model and 
we can minimize it to get the minimum DPD estimator; there is no need for kernel smoothing. 
The asymptotic properties of the minimum DPD estimators are well-studied in the literature 
and are valid under continuous models also.

Now, suppose we apply the approach of smoothed density (Basu and Lindsay, 1994) 
to find the corresponding $\textrm{MSDE}^*$ for this special case. 
Let us denote the $\textrm{MSDE}^*$ under this particular case of $\lambda = 0$ by
the minimum DPD$^*$ estimator ($\textrm{MDPDE}^*$). 
The corresponding estimating equation is given by
\begin{eqnarray}
\int (f_\theta^*)^\alpha g_n^* \widetilde{u_\theta} - \int (f_\theta^*)^{1+\alpha} \widetilde{u_\theta} &=& 0.
\end{eqnarray}
However, here we can further simplify this estimating equation by noting that
\begin{eqnarray}
\int  g_n^*(x) f_\theta^*(x)^\alpha \widetilde{u_\theta}(x)dx 
&=&  \int \left[\int W(x, y, h)dG_n(y)\right] f_\theta^*(x)^\alpha \widetilde{u_\theta}(x)dx   \nonumber\\
&=&  \int \left[\int f_\theta^*(x)^\alpha \widetilde{u_\theta}(x) W(x, y, h) dx\right]dG_n(y)  \nonumber\\
&& ~~~~~~~~~ \mbox{ (By Fubini's theorem)} \nonumber\\
&=&  \int  u_{\theta}^{\alpha*}(y) dG_n(y)  \nonumber\\
&=& \frac{1}{n} \sum_{i=1}^n  u_\theta^{\alpha*}(X_i). \nonumber
\end{eqnarray}
Here, we have used the notation 
\begin{equation}\label{EQ:u_theta_alpha_star}
u_{\theta}^{\alpha*}(y) =  \int \widetilde{u_\theta}(x) \{ f_\theta^*(x) \}^\alpha W(x,y,h) dx,
\end{equation}
in the spirit of Basu and Lindsay (1994). 
Further, taking expectation with respect to the density $f_\theta$, 
\begin{eqnarray}
E_\theta[ u_{\theta}^{\alpha*}(Y)]      
&=& \int \left[\int f_\theta^*(x)^\alpha \widetilde{u_\theta}(x) W(x, y, h) dx\right] f_\theta(y) dy  
\nonumber\\ 
&=& \int f_\theta^*(x)^\alpha \widetilde{u_\theta}(x) \left[\int W(x, y, h) f_\theta(y) dy \right] dx 
\nonumber\\
&=& \int \left(f_\theta^*(x)^\alpha \widetilde{u_\theta}(x)\right) f_\theta^*(x) dx \nonumber\\
&=& \int f_\theta^*(x)^{1+\alpha} \widetilde{u_\theta}(x).\nonumber
\end{eqnarray}
Thus, the estimating equation of the minimum DPD$^*$ estimator becomes
\begin{eqnarray}
\frac{1}{n} \sum_{i=1}^n  u_\theta^{\alpha*}(X_i)  - E_\theta [u_\theta^{\alpha*}(Y)] &=& 0, 
\label{EQ:est-eqn-MDPD_star}
\end{eqnarray}
which is clearly an unbiased estimating equation.  
Therefore the standard asymptotic results for unbiased estimating equations also hold 
for the minimum DPD$^*$ estimators. 
Under certain restrictions on the kernel function, 
the results of Section \ref{SEC:4Trans_kernel} will show that these two estimating equations 
(\ref{EQ:est-eqn-MDPD}) and (\ref{EQ:est-eqn-MDPD_star}) become the same, producing identical estimators.

\section{Influence Function of the Minimum $S^*$-Divergence Estimator}\label{SEC:MS*DE_IF}

Consider the general case of the  minimum $S^*$-divergence estimators and 
the kernel smoothed version $g^*$ of the true density $g$ defined by 
\begin{equation}
g^*(x) = \int W(x, y, h) dG(y).
\label{EQ:4smoothed_true_density}
\end{equation}
Then the minimum $S^*$-divergence estimator functional $T_{(\alpha, \lambda)}^*(G)$ 
is defined by the relation \index{Statistical Functional!for MSDE*}
$$S_{(\alpha,\lambda)}(g^*, f_{T_{(\alpha,\lambda)}^*(G)}^*) = \min_{\theta \in \Theta} 
S_{(\alpha,\lambda)}(g^*, f_\theta^*),$$
provided such a minimum exists. Then the influence function of $T_{(\alpha, \lambda)}^*$ 
can be derived from the corresponding estimating equation 
(which is the same as Equation \ref{EQ:S-div_est_eq_kernel} with $g_n^*$ replaced by $g^*$).

Let us now consider the contaminated distribution 
$G_\epsilon(x) = (1 - \epsilon) G(x) + \epsilon \wedge_y(x)$ with 
$\wedge_y$ being the degenerate distribution at the contamination point $y$. 
Suppose $g_\epsilon$ denotes the corresponding contaminated density, 
and  $g_\epsilon^*$ denotes the associated smoothed density. 
Note that $g_\epsilon^*=(1 - \epsilon) g^* + \epsilon W(x, y, h)$.
Then, by definition, the minimum $S^*$-divergence estimator functional $T(G_\epsilon)$ satisfies
\begin{equation}
\int K(\delta_\epsilon(x)) f_{\theta_\epsilon}^*(x)^{1+\alpha} \widetilde{u_{\theta_\epsilon}}(x)  dx= 0,
\label{EQ:4basu_lindsay_influence_function}
\end{equation}
where $\delta_\epsilon(x) = g_\epsilon^*(x)/f_{\theta_\epsilon}^*(x)-1$. 
Then the influence function of the $\textrm{MSDE}^*$s at the
distribution $G$ can be computed by taking a derivative of both sides
of (\ref{EQ:4basu_lindsay_influence_function}) and is presented in the following theorem. 
\index{MSDE*!Influence Function}\index{IF!of MSDE*} 

\bigskip
\begin{theorem}
	The Influence Function for the  minimum $S^*$-divergence estimator functional 
	$T_{(\alpha, \lambda)}^*$ at the distribution $G$ is given by
	\begin{eqnarray}
	IF(y; G, T_{(\alpha, \lambda)}^*) = [J_g^{*}]^{-1} N_g^{*}(y)
	\end{eqnarray}
	where
	\begin{eqnarray}
&&	N_g^{*}(y) = N_{(\alpha, \lambda)}^{*}(y; g)  \nonumber \\
	&=&  A \left[ \int (f_{\theta^g}^*(x))^B(g^*(x))^{A-1} \widetilde{u_{\theta^g}}(x) W(x,y,h)dx 
	-\int (f_{\theta^g}^*)^B(g^*)^A \widetilde{u_{\theta^g}} \right], \nonumber
	\label{EQ:4cont_N*} 
	\end{eqnarray}
	and
	\begin{eqnarray}
&&	J_g^{*} = J_{(\alpha, \lambda)}^{*}(g) \nonumber \\
	&=& A \int (f_{\theta^g}^*)^{1+\alpha} \widetilde{u_{\theta^g}}
	\widetilde{u_{\theta^g}}^T  + \int (\widetilde{i_{\theta^g}} - B \widetilde{u_{\theta^g}}
	\widetilde{u_{\theta^g}}^T )[(g^*)^A-(f_{\theta^g}^*)^A](f_{\theta^g}^*)^B, 
	~~~~~~~\label{EQ:4cont_J*} 
	\end{eqnarray}
	with ${\theta^g} = T_{(\alpha, \lambda)}^*(G)$ being the best fitting parameter under $G$ 
	and $\widetilde{i_\theta}(x) = -\nabla [ \widetilde{u_\theta}(x)]$.\hfill{$\square$}
\end{theorem}

\bigskip

\begin{corollary}\label{COR:4IF_MSDE*_model}
	When the true density $g$ belongs to the model family  
	$\{ f_\theta : \theta \in \Theta \}$,  i.e., $ g = f_\theta $ for some $\theta\in\Theta$, 
	then the Influence Function for the  Minimum $S^*$-divergence Estimator functional
	$T_{(\alpha, \lambda)}^*$ at the distribution $G = F_\theta$ becomes
	\begin{eqnarray}
	IF(y;F_\theta,T_{(\alpha, \lambda)}^*) = [J^{*}(\theta)]^{-1} \left\{ u_{\theta}^{\alpha*}(y) - 
	E_\theta [ u_{\theta}^{\alpha*}(X) ] \right\},
	\end{eqnarray}
	where
	\begin{eqnarray}\label{EQ:4cont_J*_model} 
	J^* = J_{(\alpha, \lambda)}^{*}(f_\theta)  &=& E_\theta [ u_{\theta}^{2\alpha*}(X) ],
	\end{eqnarray}
	with
	\begin{eqnarray}\label{EQ:u_theta_2alpha_star}
	u_{\theta}^{2\alpha*}(y) &=&  \int \widetilde{u_\theta}(x) 
	\widetilde{u_\theta}(x)^T \{ f_\theta^*(x) \}^\alpha W(x,y,h) dx. 
	\end{eqnarray}
	Note that the Influence Function at the model distribution turns out to be independent 
	of the parameter $\lambda$, just as it was for the minimum $S$-divergence  estimator  
	(see Equation \ref{EQ:S_div_IF_model}).\hfill{$\square$}
\end{corollary}

\begin{remark}
	In is interesting to note that, the matrix $J^*$ defined in (\ref{EQ:4cont_J*_model}) can also be written as 
	\begin{eqnarray}
	J^* = E_\theta [ u_{\theta}^{2\alpha*}(X) ] = E_\theta [ - \nabla u_{\theta}^{\alpha*}(X) ].
	\label{EQ:4cont_J*_model2}
	\end{eqnarray}
	To see this, note that 
	$$
	\nabla u_{\theta}^{\alpha*}(X) = \int \nabla \widetilde{u_\theta}(x)\{ f_\theta^*(x) \}^\alpha W(x,X,h) dx
	+ \alpha u_{\theta}^{2\alpha*}(X),
	$$
	and 
	\begin{eqnarray}
	&& E_\theta \left[\int \nabla \widetilde{u_\theta}(x)\{ f_\theta^*(x) \}^\alpha W(x,X,h) dx\right]\nonumber\\
	&=& \int \left[\int \nabla \widetilde{u_\theta}(x)\{ f_\theta^*(x) \}^\alpha W(x,y,h) dx\right]f_\theta(y) dy
	\nonumber\\
	&=& \int \nabla \widetilde{u_\theta}(x)\{ f_\theta^*(x) \}^\alpha \left[\int f_\theta(y) W(x,y,h) dy\right]dx
	\nonumber\\
	&=& \int \nabla \widetilde{u_\theta}(x)\{ f_\theta^*(x) \}^{1+\alpha} dx
	\nonumber\\
	&=& -(1+\alpha) \int \widetilde{u_\theta}(x) \widetilde{u_\theta}(x)^T \{ f_\theta^*(x) \}^{1+\alpha} dx,
	\end{eqnarray}
	where the last step follows using integration by parts.
	A similar argument shows that
	\begin{eqnarray}
	J^* = E_\theta [ u_{\theta}^{2\alpha*}(X) ] 
	= \int \widetilde{u_\theta}(x) \widetilde{u_\theta}(x)^T \{ f_\theta^*(x) \}^{1+\alpha} dx.
	\label{EQ:4cont_J*_model3}
	\end{eqnarray}
	Combining all these, we get 
	$$
	E_\theta\left[\nabla u_{\theta}^{\alpha*}(X)\right] 
	= E_\theta\left[ \int \nabla \widetilde{u_\theta}(x)\{ f_\theta^*(x) \}^\alpha W(x,X,h) dx\right] 
	+ \alpha E_\theta\left[u_{\theta}^{2\alpha*}(X)\right],
	$$
	or, 
	$$
	E_\theta\left[\nabla u_{\theta}^{\alpha*}(X)\right]  =  -(1+\alpha) J^* + \alpha J^* = - J^*.
	$$
	This proves (\ref{EQ:4cont_J*_model2}).
	\hfill{$\square$}
\end{remark}


\section{Asymptotic Distribution of the Minimum $S^*$-Divergence Estimator}\label{SEC:MS*DE_asymptotic}

Now, we will derive the asymptotic properties of the minimum $S^*$-divergence estimator 
and prove its asymptotic equivalence with the minimum $S$-divergence estimators under some assumptions 
on the kernel used. Consider the parametric set-up of the continuous model 
$\{f_\theta : \theta \in \Theta\}$ with corresponding smoothed model 
$f_\theta^*$ as defined in equation (\ref{EQ:3smoothed_model}). 
We have $n$ independent and identically distributed observations 
$X_1, \ldots, X_n$ from the true density $g$ and the corresponding smoothed density $g^*$ is as defined 
in equation (\ref{EQ:4smoothed_true_density}). Along with all the notations of the previous sections, we 
will further assume some more conditions to give a rigorous proof of 
the consistency and asymptotic normality of the minimum $S^*$-divergence estimator.  
The assumptions will be based on the following general definitions 
about the properties of the model densities.

\begin{definition} \label{DFN:1identifiable}
	The parametric model ${\cal F}$ is said to be identifiable, 
	if for any  $\theta_1, \theta_2 \in \Theta$, $\theta_1
	\neq \theta_2$ implies 
	$f_{\theta_1}(x) \neq f_{\theta_2}(x)$ on a set of
	positive dominating measure, where $f_\theta$ represents the density
	function of $F_\theta$.
\end{definition}

\begin{definition} \label{DFN:1compatible_with_family}
	Let the parametric densities  $\{f_\theta: \theta \in \Theta\}$
	have common support $K^*$, independent of $\theta$. Then the true density $g$ is said to be compatible
	with the family $\{f_\theta\}$ if $K^*$ is also the support of $g$.
\end{definition}

\begin{definition} \label{DFN:1smooth_kernel_integrated_family} 
	The kernel integrated family of distributions  $f_\theta^*(x)$ defined in equation
	(\ref{EQ:3smoothed_model}) is called smooth if the conditions of Lehmann (1983, p. 409, p. 429) 
	are satisfied with $f_\theta^*(x)$ is place of  $f_\theta(x)$. 
\end{definition}

Then the assumptions required for proving the asymptotic results are as follows:

\renewcommand{\theenumi}{(SB\arabic{enumi})}
\renewcommand\labelenumi{(SB\arabic{enumi})}
\begin{enumerate}
	\item\label{ASM:4asSB1} The model family $\mathcal{F}$ is identifiable 
	in the sense of Definition \ref{DFN:1identifiable}.
	
	\item\label{ASM:4asSB2} The probability density functions $f_\theta$ of the model distribution have common 
	support so that the set $\mathcal{X} = \{ x : f_{\theta}(x) >0 \}$ is independent of $\theta$.
	Also the true distribution $g$ is compatible with the model family in the sense of 
	Definition  \ref{DFN:1compatible_with_family}.
	
	\item\label{ASM:4asSB3} The kernel integrated family of distributions  $f_\theta^*(x)$ is smooth
	in the sense of Definition \ref{DFN:1smooth_kernel_integrated_family}.
	
	\item\label{ASM:4asSB4} The matrix $J_g^{*}(\theta^g)$ is positive definite.
	
	\item\label{ASM:4asSB5} The quantities  
	$\int (g_n^*)^{1/2}(x) (f_{\theta}^*(x))^{\alpha} |\widetilde{ u_{j\theta}}(x)|dx$, \\~~~~
	$\int (g_n^*)^{1/2}(x) (f_{\theta}^*(x))^{\alpha}|\widetilde{ u_{j\theta}}(x)| |\widetilde{ u_{k\theta}}(x)|dx$  
	and		 $\int (g_n^*)^{1/2}(x) (f_{\theta}^*(x))^{\alpha} |\widetilde{ u_{jk\theta}}(x)|dx$ \\
	are bounded for all $j,k$ and for all $\theta \in \omega$, an open neighborhood of the best fitting parameter $\theta^g$.

	\item\label{ASM:4asSB6} For almost all $x$, there exists functions $M_{jkl}(x)$, $M_{jk,l}(x)$, 
	$M_{j,k,l}(x)$ that dominate, in absolute value, \begin{center}
		$\left[(f_{\theta}^*(x))^{\alpha} \widetilde{ u_{jkl\theta}}(x)\right]$, ~
		$\left[(f_{\theta}^*(x))^{\alpha} \widetilde{ u_{jk\theta}}(x)\widetilde{ u_{l\theta}}(x)\right]$ ~
		and $\left[(f_{\theta}^*(x))^{\alpha} \widetilde{ u_{j\theta}}(x) 
		\widetilde{ u_{k\theta}}(x) \widetilde{ u_{l\theta}}(x)\right]$ 
	\end{center}
	respectively for all $j, k, l$ and are uniformly bounded in expectation with respect to $g^*$ 
	and $f_{\theta}^*$  for all $\theta \in \omega$.
	
	\item\label{ASM:4asSB7} The function $\left( \frac{g^*(x)}{f_{\theta}^*(x)} \right)^{A-1}$ is uniformly 
	bounded (by, say, $C$) for all $\theta \in \omega$.
\end{enumerate}
\renewcommand{\theenumi}{\arabic{enumi}}
\renewcommand\labelenumi{\arabic{enumi}}

To prove the desired asymptotic results for the minimum $S^*$-divergence estimator, 
we will, from now on, assume that the conditions \ref{ASM:4asSB1}--\ref{ASM:4asSB7} hold. 
Also, as in the case of 
asymptotic distribution of the MSDE under discrete models, here also we will first prove some important Lemmas. 
Let us define $\eta_n(x) = \sqrt{n} (\sqrt{\delta_n^*} -\sqrt{\delta_g^*})^2$ with
\begin{eqnarray}
\delta_g^* &=& \frac{g^*(x)}{f_{\theta}^*(x)} - 1. \nonumber
\end{eqnarray}
\begin{lemma}\label{LEM:4cont_asymp_lem1}
	Provided it exists, $Var_g(g_n^*(x)) = \frac{1}{n} \nu(x)$, 
	where $\nu(x)$ is defined as
	\begin{equation}
	\nu(x) = \int W^2(x,y,h)g(y)dy - [g^*(x)]^2.
	\end{equation}
	\hfill{$\square$}
\end{lemma}
For the proof of this Lemma see Basu and Lindsay (1994). 
Next we need to assume that the boundedness of the kernel function $W(x,y,h)$. 
So, from now on, let
$$
W(x, y, h) \le M(h) < \infty,
$$
where $M(h)$ depends on $h$, but not on $x$ or $y$. Then, we have
\begin{eqnarray}
\nu(x) &\le&   \int W^2(x,y,h)g(y)dy  \nonumber \\ \nonumber  
&\le&	M(h) \int W(x,y,h)g(y)dy  \\ 
&=& M(h)g^*(x).  \label{EQ:4cont_nu}
\end{eqnarray}
The next three Lemmas follow from Basu and Lindsay (1994) and are stated here without proof.

\begin{lemma}
	\label{LEM:4cont_asymp_lem2}
	Under the above mentioned set-up and notations, 
	$$
	n^{1/4}\left[(g_n^*(x))^{1/2} - g^*(x)^{1/2}\right] \rightarrow 0,
	$$
	with probability $1$, provided $ \nu(x) < \infty$.
	\hfill{$\square$}
\end{lemma} 
\begin{lemma}\label{LEM:4cont_asymp_lem3}
	For any $k \in [0,2]$, we have
	\begin{enumerate}
		\item $E[\eta_n^*(x)^k] \le n^{\frac{k}{2}}E[|\delta_n^*(x) - \delta_g^*(x)|]^k \le 
		\left[ \frac{\nu(x)}{(f_{\theta}^*(x))^2}\right]^{\frac{k}{2}}$. \\
		\item   $E[|\delta_n^*(x) - \delta_g^*(x)|] \le \frac{\sqrt{\nu(x)}}{f_{\theta}^*(x)}$.
	\end{enumerate}
	\hfill{$\square$}
\end{lemma}

\begin{lemma}\label{LEM:4cont_asymp_lem4}
	$E[\{\eta_n^*(x)\}^k] \rightarrow 0$, as $n \rightarrow \infty$,  for $k \in [0,2).$
	\hfill{$\square$}
\end{lemma}
\medskip
\noindent
Next, let us define, $$a_n^*(x) = K(\delta_n^*(x)) - K(\delta_g^*(x)),$$ 
$$b_n^*(x) = (\delta_n^*(x)-\delta_g^*(x))K'(\delta_g^*(x))$$ 
$$\mbox{and~~~~}\tau_n^*(x) = \sqrt n |a_n^*(x) - b_n^*(x)|.$$ 
We will need the limiting distributions of
$$S_{1n}^* = \sqrt n \sum_x a_n^*(x)(f_{\theta}^*(x))^{1+\alpha} \widetilde{ u_{\theta}}(x),~~~~ 
\mbox{and~~~~}S_{2n}^* = \sqrt n \sum_x b_n^*(x)(f_{\theta}^*(x))^{1+\alpha}  \widetilde{ u_{\theta}}(x).$$
\begin{lemma}\label{LEM:4cont_asymp_lem5}
	Assume condition \ref{ASM:4asSB5}. Then 
\begin{eqnarray}
E|S_{1n}-S_{2n}| \rightarrow 0, ~~~~ &\mbox{as}& ~~~~~ n \rightarrow \infty,\nonumber\\ 
\mbox{and hence ~~~~~~~~~~~} S_{1n}-S_{2n} \displaystyle\mathop{\rightarrow}^\mathcal{P} 0, ~~~~ 
	&\mbox{as}& ~~~~~ n \rightarrow \infty.\nonumber
\end{eqnarray}
\end{lemma}
\noindent\textbf{Proof: }
By Lemma 25 of Lindsay (1994),
there exists some positive constant $\beta$ such that 
$$
\tau_n^*(x) \le \beta \sqrt{n} (\sqrt{\delta_n^*} - \sqrt{\delta_g^*})^2 = \beta \eta_n^*(x).
$$ 
Also, by Lemma \ref{LEM:4cont_asymp_lem3} and equation (\ref{EQ:4cont_nu}), 
$$ E[\tau_n^*(x)] \le \beta \frac{\nu^{1/2}(x)}{f_{\theta}^*(x)}
\le \beta \frac{M^{1/2}(h)(g^*(x))^{1/2}}{f_{\theta}^*(x)}.$$
By Lemma \ref{LEM:4cont_asymp_lem4}, 
$$E[\tau_n^*(x)] = \beta E[\eta_n^*(x)] \rightarrow 0~~~\mbox{ as }~~n \rightarrow \infty.$$ 
Thus we get, 
\begin{eqnarray}
E|S_{1n}^*-S_{2n}^*| &\le& \sum_x E[\tau_n^*(x)] (f_{\theta}^*(x))^{1+\alpha}|\widetilde{ u_{\theta}}(x)| 
\nonumber \\
&\le&  \beta M^{1/2}(h) \sum_x (g^*(x))^{1/2} (f_{\theta}^*(x))^{\alpha}|\widetilde{ u_{\theta}}(x)| \nonumber \\
& < &   \infty, ~~~~~~ \mbox{ (by assumption SB5).} \nonumber 
\end{eqnarray}
So, by the Dominated Convergence Theorem (DCT),  
$$E|S_{1n}^* - S_{2n}^*| \rightarrow 0 ~~~\mbox{ as }~~ n \rightarrow \infty.$$
Hence, by Markov inequality, 
$$S_{1n}^*-S_{2n}^* \displaystyle\mathop{\rightarrow}^\mathcal{P} 0 ~~~\mbox{ as }~~ n \rightarrow \infty.$$
\hfill{$\square$}
\begin{lemma}\label{LEM:4cont_asymp_lem6}
	Suppose, 
	\begin{equation}\label{EQ:4cont_V*} 
	V_g^* = V_{(\alpha, \lambda)}^*(g) =  Var \left[ \int W(x,X,h)K'(\delta_g^{g*}(x)) (f_{\theta^g}^*(x))^\alpha 
	\widetilde{ u_{\theta^g}}(x) dx \right]
	\end{equation}
	is finite. Then under the true distribution $g$, 
	$$
	S_{1n}^* \displaystyle\mathop{\rightarrow}^\mathcal{D} N(0, V_g^*)
	$$
\end{lemma}
\noindent\textbf{Proof: }
Note that, by Lemma \ref{LEM:4cont_asymp_lem5}, the asymptotic distribution 
of $S_{1n}^*$ and $S_{2n}^*$ are the same. Now, we have
\begin{eqnarray}
S_{2n}^* &=& \sqrt n \sum_x (\delta_n(x) - \delta_g(x)) K'(\delta_g(x)) (f_{\theta}^*(x))^{1+\alpha}
\widetilde{ u_{\theta}}(x)\nonumber \\
&=&  \sqrt n \sum_x (g_n^*(x) - g^*(x)) K'(\delta_g(x)) (f_{\theta}^*(x))^\alpha
\widetilde{ u_{\theta}}(x) \nonumber \\
&=&  \sqrt n \left( \frac{1}{n} \sum_{i=1}^n \int \left[ W(x,X_i,h) - E\{W(x,X_i,h)\} \right]
K'(\delta_g(x))(f_{\theta}^*(x))^\alpha \widetilde{ u_{\theta}}(x) dx \right) \nonumber \\
&\displaystyle\mathop{\rightarrow}^\mathcal{D} & N( 0, V) ~~~~~ \mbox{(by the central limit theorem).} \nonumber 
\end{eqnarray}
This completes the proof of the Lemma. 
\hfill{$\square$}

\bigskip
Now, using the above lemmas,  we can prove the main theorem of this section regarding the 
consistency and asymptotic normality of the Minimum $S^*$-divergence  Estimator.\\
\begin{theorem}\label{THM:4cont_asymp}
	Under Assumptions \ref{ASM:4asSB1}--\ref{ASM:4asSB7}, 
	there exists a consistent sequence $\theta_n^*$ of roots to the 
	Minimum $S^*$-divergence estimating equation (\ref{EQ:S-div_est_eq_kernel}).\\
	Also, the asymptotic distribution of $\sqrt n (\theta_n^* - \theta^g)$ is $p-$dimensional 
	normal with mean $0$ and variance $[J_g^{*}]^{-1}V_g^* [J^{*}_{g}]^{-1}$,
	where $J_g^*$ and $V_g^*$ are as defined in equation (\ref{EQ:4cont_J*}) and (\ref{EQ:4cont_V*}).
\end{theorem}
\noindent\textbf{Proof: }
The proof is similar to that of the minimum $S$-divergence estimators under discrete 
models (Ghosh, 2014) with slight modifications based on Lemmas 
\ref{LEM:4cont_asymp_lem1} to \ref{LEM:4cont_asymp_lem6}.
\hfill{$\square$}

\bigskip
\begin{corollary}\label{COR:4cont_asymp_model1}
	When the true density $g$ belongs to the model family  $\{ f_\theta : \theta \in \Theta \}$, i.e.,
	$ g = f_\theta $, the asymptotic distribution of $\sqrt n (\theta_n^* - \theta)$  is 
	Normal with mean zero and variance-covariance matrix $[J^{*}]^{-1}V^* [J^{*}]^{-1}$, 
	where $V^* = V_{(\alpha, \lambda)}^*(f_\theta) = V_\theta[u_\theta^{\alpha*}(X)]$ 
	and $J^*$ is defined in equation (\ref{EQ:4cont_J*_model}).
	Note that this asymptotic distribution turns out to be independent of the parameter $\lambda$.  
	(This independence is in-line with that of the minimum $S$-divergence estimator for discrete models 
	derived in Ghosh, 2014).\hfill{$\square$}
\end{corollary}

\bigskip 
\begin{remark}\label{REM:4asymp_var_H0}
	It is to be noted that the result of Theorem \ref{THM:4cont_asymp} is a fixed bandwidth ($h$) results. 
	Further, if we let $h \rightarrow 0$, then we can easily show that 
	under appropriate regularity conditions the asymptotic variance matrix 
	$[J^{*}]^{-1}V^* [J^{*}]^{-1}$ at the model density converges element-wise to matrix 
	\begin{eqnarray}
	\left( E_\theta[u_\theta(X)u_\theta(X)^T f_\theta^\alpha(X)] \right)^{-1} 
	V_\theta[u_\theta(X)f_\theta^\alpha(X)] 
	\left( E_\theta[u_\theta(X)u_\theta(X)^T f_\theta^\alpha(X)] \right)^{-1}. 
	\label{EQ:4cont_asymp_var_model}
	\end{eqnarray}
	Interestingly this matrix (\ref{EQ:4cont_asymp_var_model}) is the same as the 
	asymptotic variance of the minimum $S$-divergence estimators at the discrete model,
	as obtained in Ghosh (2014).

	In the same spirit, we also expect that under the additional condition $h \rightarrow 0$
	along with the assumptions of Theorem \ref{THM:4cont_asymp}, 
	the asymptotic distribution of $\sqrt(n) (\theta_n^* - \theta^g)$
	will be $p$-variate normal with mean vector $0$ and variance given by expression (\ref{EQ:4cont_asymp_var_model}), 
	which corresponds to the discrete case, although we do not have a fully rigorous proof at this point.

	As an illustration of this expectation, we compute the asymptotic variance 
	$[J^{*}]^{-1}V^* [J^{*}]^{-1}$ of the minimum $S^*$-divergence estimators at the model 
	for $N(\theta, \sigma^2)$ model distribution with several values of the fixed bandwidth $h$ converging to zero and
	compare these values with the asymptotic variance given in expression (\ref{EQ:4cont_asymp_var_model}).
	We consider the Gaussian kernel with $W(x, y, h)$ being the $N(y, h^2)$ density at $x$
	and different values of $\alpha\geq 0$. Then, a simple calculation shows that 
	the smoothed density $f_\theta^* (x)$ is the normal $N(\theta, \sigma^2+h^2)$ density
	and the corresponding asymptotic variance at the model given in Corollary \ref{COR:4cont_asymp_model1}
	simplifies to 
	$$
	[J^{*}]^{-1}V^* [J^{*}]^{-1} = \left[\frac{(1+\alpha)^2(\sigma^2+h^2)^2}{
		((1+\alpha)h^2 + \sigma^2)((1+\alpha)h^2 + (1+2\alpha)\sigma^2)}\right]^{3/2}\sigma^2
	= \zeta_{\alpha, h}\sigma^2.
	$$
	Interestingly, in this case, we have $\zeta_{0, h} = 1$ so that 
	the value of $[J^{*}]^{-1}V^* [J^{*}]^{-1}$ at $\alpha=0$ further simplifies to $\sigma^2$ 
	which is independent of $h$ and also the same as the corresponding value of 
	expression (\ref{EQ:4cont_asymp_var_model}); we will examine this special 
	property towards the end of this paper in Section \ref{SEC:4Trans_kernel}.
	For all other $\alpha >0$, the value of the expression (\ref{EQ:4cont_asymp_var_model}) 
	equals $\zeta_\alpha\sigma^2$ with $\zeta_\alpha = (1+\alpha)^3(1+2\alpha)^{-3/2}$
	and clearly 
	$$
	\zeta_{\alpha, h} \rightarrow \zeta_{\alpha}, ~~~~~~ \mbox{ as }~~h \rightarrow 0.
	$$ 
	The expression $\zeta_{\alpha} \sigma^2$ is the asymptotic variance of the MDPDE of $\mu$.
	This gives some substance to our description in the previous paragraphs relating to 
	the asymptotic variance of the MSDE$^*$ under the condition $h \rightarrow 0$.  
\hfill{$\square$}
\end{remark}


\section{Simulation Studies: Normal Model}\label{SEC:4MS*DE_simulation}

We will now explore the performance of the minimum $S^*$-divergence estimator 
through a detailed simulation study. For simplicity, we will consider the model density $f_\theta$ 
to be the normal density with unknown mean $\mu$ and unknown variance $\sigma^2$ so that 
the parameter of interest is $\theta=(\mu, ~ \sigma)$ and the parameter space is
$\Theta=\mathbb{R}\times[0,\infty)$. We will simulate sample data of size $n=50$ 
from a normal distribution with  mean $0$ and variance $3$ and compute the minimum $S^*$-divergence
estimators based on the sample drawn. We compute the empirical bias and MSE for the minimum 
$S^*$-divergence estimator over $1000$ replications. For computation of the 
minimum $S^*$-divergence estimator under the Basu-Lindsay approach, we will 
use the Gaussian kernel and, following  the normal reference rule (Scott, 2001), use the bandwidth 
$$
h_n = 1.06 \sigma n^{-1/5},
$$
with $\sigma$ being the standard deviation. In this expression of the bandwidth, 
we will replace $\sigma$ by its robust estimate $\sigma_0$ defined in Equation (\ref{EQ:sigma_0_est}).
Also, we need to form the smoothed model density as described in previous sections.
Here our model is the $N(\mu, \sigma^2)$ density
and we choose the Gaussian kernel with $W(x, y, h)$ being the $N(y, h^2)$ density at $x$.
We have already observed that the smoothed density $f_\theta^* (x)$ is 
the normal $N(\mu, \sigma^2+h^2)$ density in this case.

\begin{table}[h]
	\centering 
	\caption{MSEs of the location parameter $\mu$ without any contamination }
	\begin{tabular}{r| lllllll} \hline\noalign{\smallskip}
		$\lambda$	&	$\alpha = $ 0	&	$\alpha = $ 0.1	&	$\alpha = $ 0.3	&	$\alpha = $ 0.4	&
		$\alpha = $ 0.5	&		$\alpha = $ 0.8	&	$\alpha = $ 1	\\ 
		\noalign{\smallskip}\hline\noalign{\smallskip}
		$ -1$	&	0.22	&	0.17	&	0.17	&	0.18	&	0.19	&	0.21	&	0.22	\\
		$-0.7$	&	0.16	&	0.17	&	0.17	&	0.18	&	0.19	&	0.21	&	0.22	\\
		$-0.5$	&	0.16	&	0.17	&	0.17	&	0.18	&	0.19	&	0.21	&	0.22	\\
		$-0.3$	&	0.16	&	0.17	&	0.17	&	0.18	&	0.19	&	0.21	&	0.22	\\
		0	&	0.16	&	0.16	&	0.17	&	0.18	&	0.19	&	0.21	&	0.22	\\
		0.5	&	0.16	&	0.16	&	0.17	&	0.18	&	0.19	&	0.21	&	0.22	\\
		1	&	0.17	&	0.17	&	0.17	&	0.18	&	0.18	&	0.21	&	0.22	\\
		1.5	&	0.17	&	0.17	&	0.17	&	0.18	&	0.18	&	0.21	&	0.22	\\
		2	&	0.17	&	0.17	&	0.17	&	0.18	&	0.18	&	0.21	&	0.22	\\
	\noalign{\smallskip}	\hline
	\end{tabular}
	\label{TAB:mu_WO}
\end{table}

\begin{table}[h]
	\centering 
	\caption{MSEs of the scale parameter $\sigma$ without any contamination }
	\begin{tabular}{r| lllllll} \hline\noalign{\smallskip}
		$\lambda$	&	$\alpha = $ 0	&	$\alpha = $ 0.1	&	$\alpha = $ 0.3	&	$\alpha = $ 0.4	&	$\alpha = $ 0.5	&	
		$\alpha = $ 0.8	&	$\alpha = $ 1	\\ 
		\noalign{\smallskip} \hline\noalign{\smallskip}
		$-1$	&	0.16	&	0.11	&	0.11	&	0.11	&	0.12	&	0.13	&	0.14	\\
		$-0.7$	&	0.11	&	0.10	&	0.11	&	0.11	&	0.12	&	0.13	&	0.14	\\
		$-0.5$	&	0.10	&	0.10	&	0.10	&	0.11	&	0.12	&	0.13	&	0.14	\\
		$-0.3$	&	0.10	&	0.10	&	0.10	&	0.11	&	0.11	&	0.13	&	0.14	\\
		0	&	0.09	&	0.09	&	0.10	&	0.11	&	0.11	&	0.13	&	0.14	\\
		0.5	&	0.10	&	0.09	&	0.10	&	0.10	&	0.11	&	0.13	&	0.14	\\
		1	&	0.10	&	0.10	&	0.10	&	0.10	&	0.11	&	0.13	&	0.14	\\
		1.5	&	0.12	&	0.11	&	0.10	&	0.10	&	0.11	&	0.13	&	0.14	\\
		2	&	0.13	&	0.12	&	0.11	&	0.10	&	0.10	&	0.13	&	0.14	\\
		\noalign{\smallskip}		\hline
	\end{tabular}
	\label{TAB:sigma_WO}
\end{table}

First we will consider the case of pure data without any contamination. So, we generate samples 
from pure $N(0,3)$ density; the computed MSE of the mean and variance parameter are reported in 
Tables \ref{TAB:mu_WO} and \ref{TAB:sigma_WO}. Clearly, as $\alpha$ increases from $0$ to $1$, 
the MSE increases slightly and so there is a loss in efficiency of the minimum $S^*$-divergence estimators
with increasing $\alpha$. However, this loss in efficiency is quite small and does not limit 
its application in pure data.

The main reason of using the minimum $S^*$-divergence estimator is that 
it gives the option of several robust estimators at different but high levels of efficiency.
From the influence function analysis of Section \ref{SEC:MS*DE_IF}, it is suggested that the
robustness of these estimators should increase with increasing $\alpha$.
To study the robustness under several types of contaminations, we next consider data from the 
$N(0,3)$ distribution but with $5\%$, $10\%$ or $20\%$ contamination from different distributions as follows:
\begin{enumerate}
	\item[(i)] mean shifted density -- $N(15,3)$
	\item[(ii)] larger variance -- $N(0,10)$
	\item[(iii)] smaller variance -- $N(0,1)$
	\item[(iii)] symmetric contamination with heavy tail -- $t$-density with $1$ degrees of freedom.
	\item[(iv)] Non-symmetric contamination with heavy tail -- $\chi^2$-density with $10$ degrees of freedom.
\end{enumerate}

\begin{figure}[!th]
	\centering
	\includegraphics[width=0.3\textwidth]{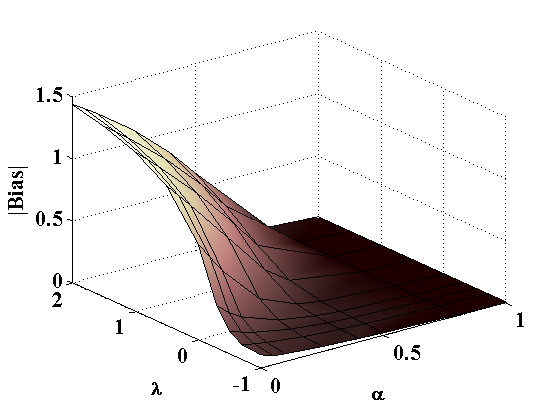}
	\includegraphics[width=0.3\textwidth]{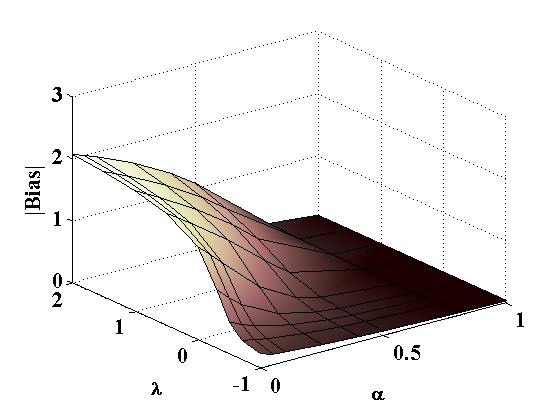}
	\includegraphics[width=0.3\textwidth]{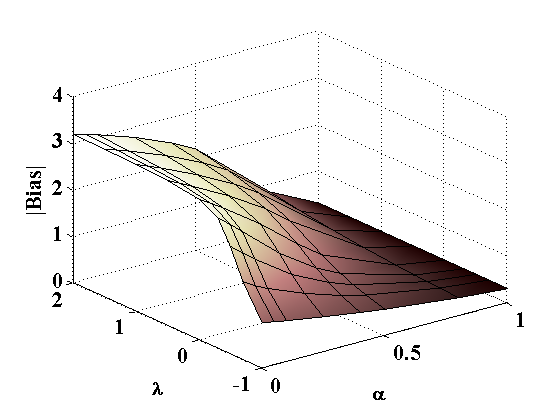}
	\\
	\subfloat[$\epsilon=5\%$]{
		\includegraphics[width=0.3\textwidth]{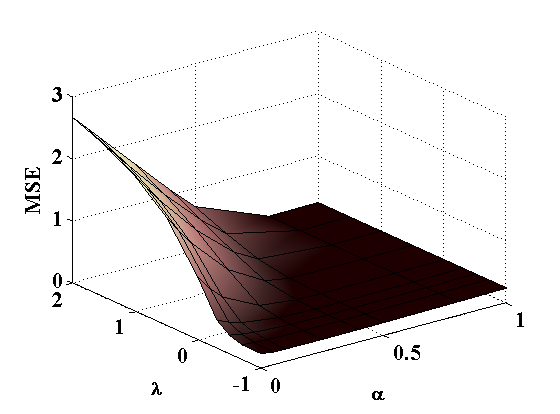}
		\label{FIG:mu_1_5}}
	\subfloat[$\epsilon=10\%$]{
		\includegraphics[width=0.3\textwidth]{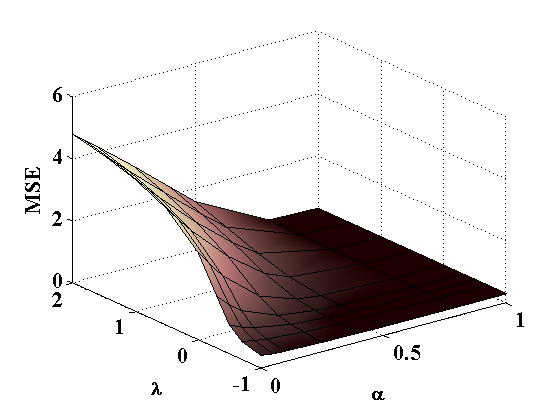}
		\label{FIG:mu_1_10}}
	\subfloat[$\epsilon=20\%$]{
		\includegraphics[width=0.3\textwidth]{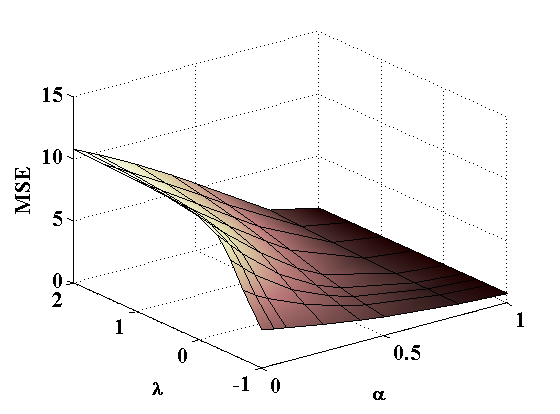}
		\label{FIG:mu_1_20}}
	\caption{Bias and MSE of the estimator of location parameter $\mu$ with different contamination proportion $\epsilon$ for case (i)}
	\label{FIG:mu_1}
\end{figure}

\begin{figure}[!th]
	\centering
	\includegraphics[width=0.3\textwidth]{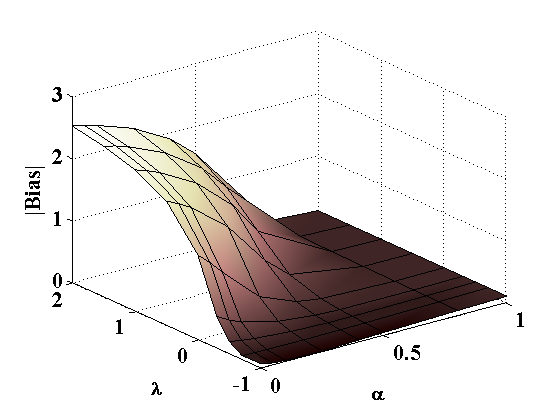}
	\includegraphics[width=0.3\textwidth]{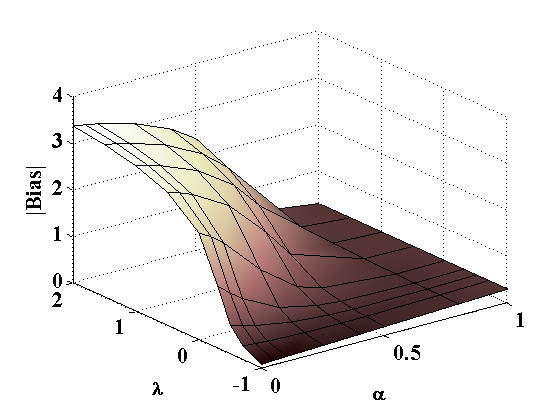}
	\includegraphics[width=0.3\textwidth]{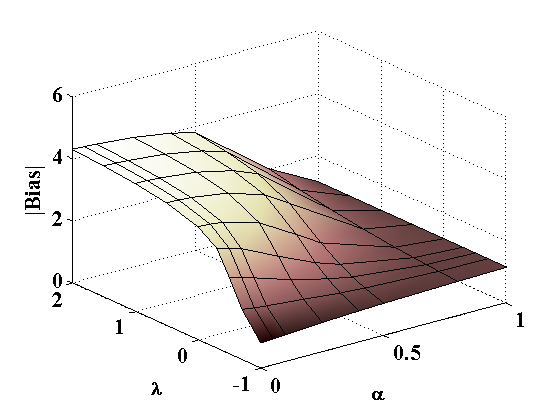}
	\\
	\subfloat[$\epsilon=5\%$]{
		\includegraphics[width=0.3\textwidth]{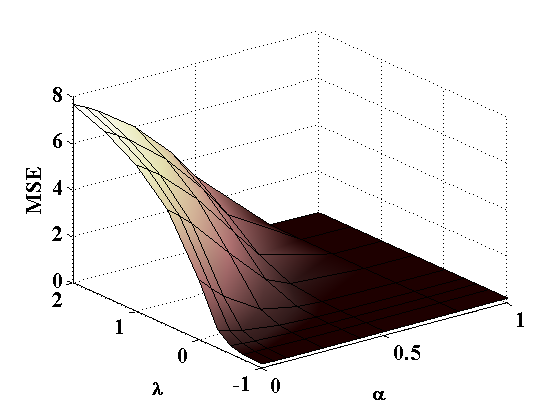}
		\label{FIG:sd_1_5}}
	\subfloat[$\epsilon=10\%$]{
		\includegraphics[width=0.3\textwidth]{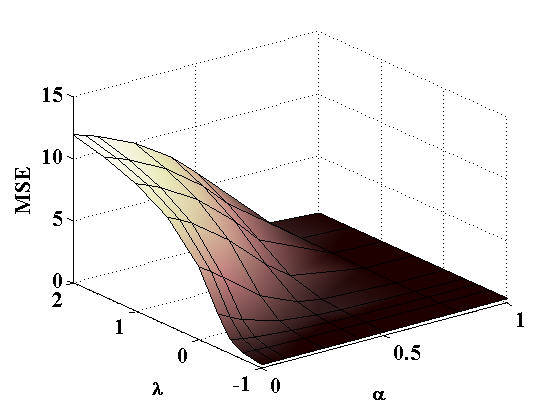}
		\label{FIG:sd_1_10}}
	\subfloat[$\epsilon=20\%$]{
		\includegraphics[width=0.3\textwidth]{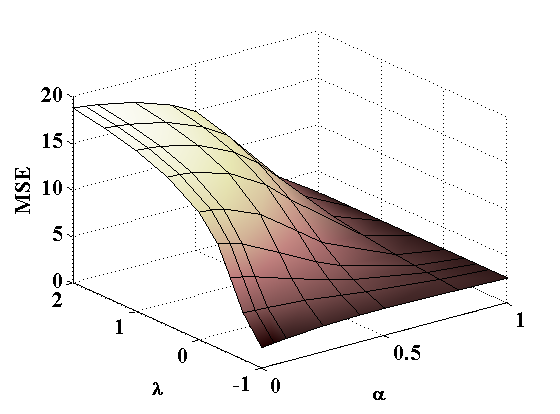}
		\label{FIG:sd_1_20}}
	\caption{Bias and MSE of the  estimator of scale parameter $\sigma$ with different contamination proportion $\epsilon$ for case (i)}
	\label{FIG:sd_1}
\end{figure}

\begin{figure}[!th]
	\centering
	\includegraphics[width=0.3\textwidth]{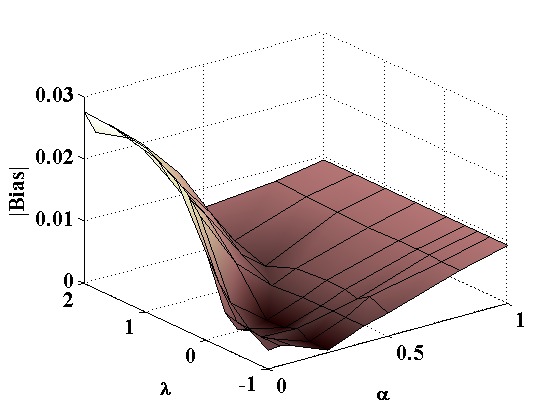}
	\includegraphics[width=0.3\textwidth]{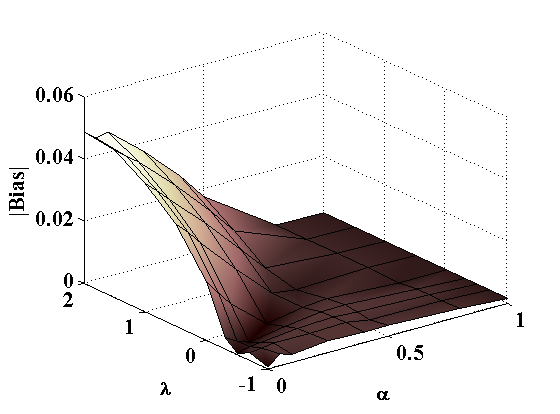}
	\includegraphics[width=0.3\textwidth]{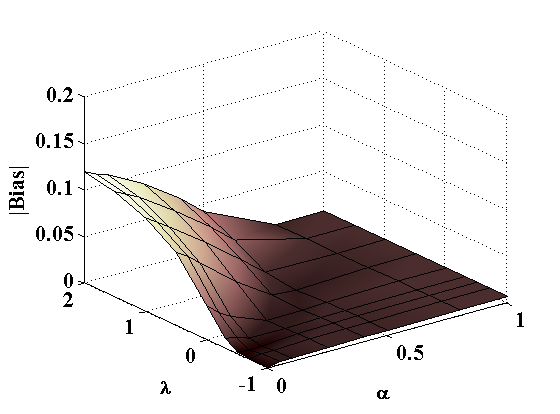}
	\\
	\subfloat[$\epsilon=5\%$]{
		\includegraphics[width=0.3\textwidth]{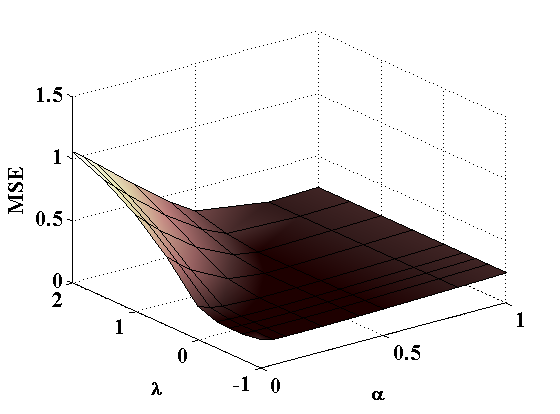}
		\label{FIG:mu_2_5}}
	\subfloat[$\epsilon=10\%$]{
		\includegraphics[width=0.3\textwidth]{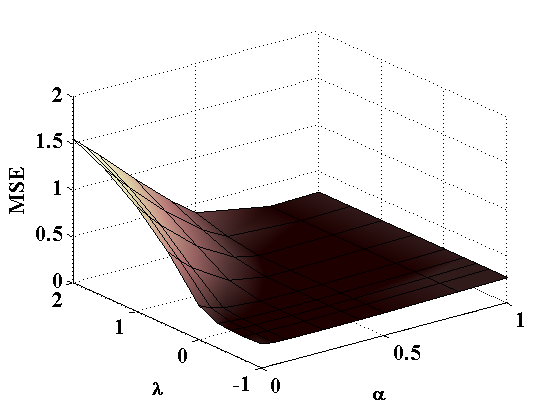}
		\label{FIG:mu_2_10}}
	\subfloat[$\epsilon=20\%$]{
		\includegraphics[width=0.3\textwidth]{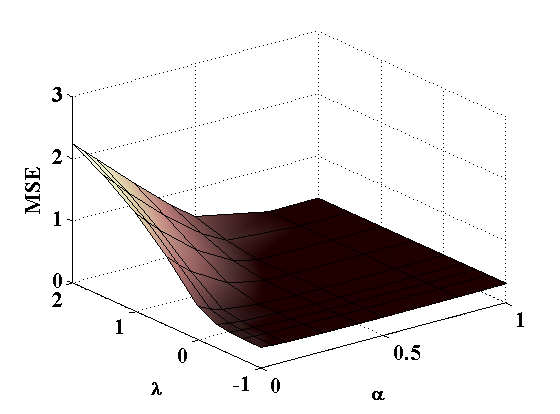}
		\label{FIG:mu_2_20}}
	\caption{Bias and MSE of the  estimator of location parameter $\mu$ with different contamination proportion $\epsilon$ for case (ii)}
	\label{FIG:mu_2}
\end{figure}

\begin{figure}[!th]
	\centering
	\includegraphics[width=0.3\textwidth]{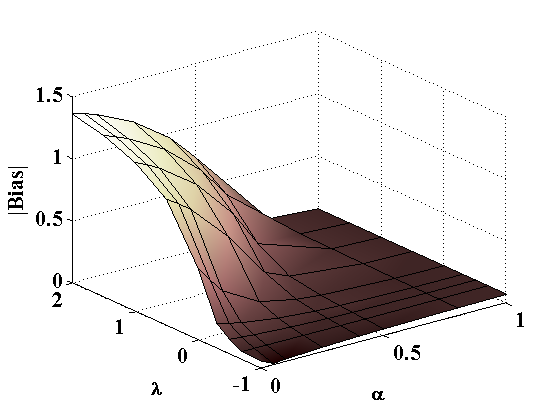}
	\includegraphics[width=0.3\textwidth]{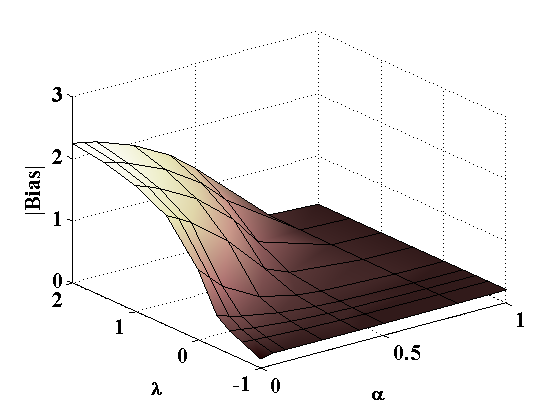}
	\includegraphics[width=0.3\textwidth]{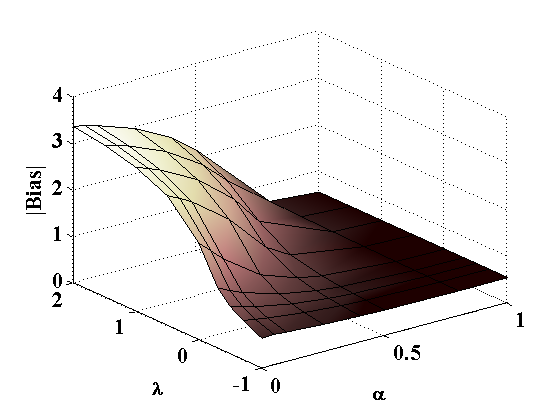}
	\\
	\subfloat[$\epsilon=5\%$]{
		\includegraphics[width=0.3\textwidth]{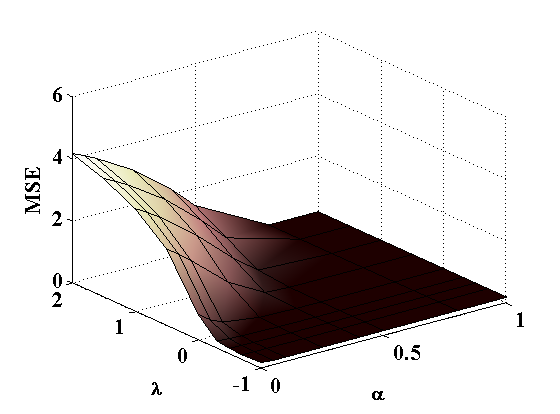}
		\label{FIG:sd_2_5}}
	\subfloat[$\epsilon=10\%$]{
		\includegraphics[width=0.3\textwidth]{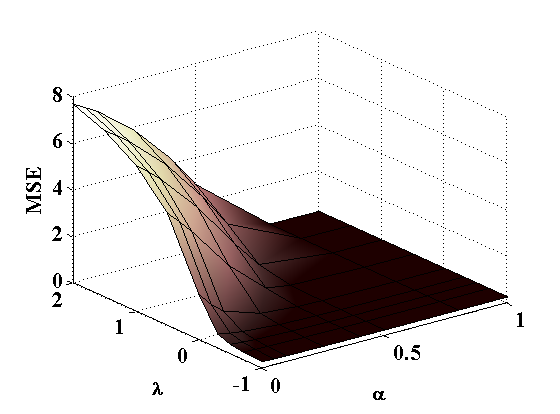}
		\label{FIG:sd_2_10}}
	\subfloat[$\epsilon=20\%$]{
		\includegraphics[width=0.3\textwidth]{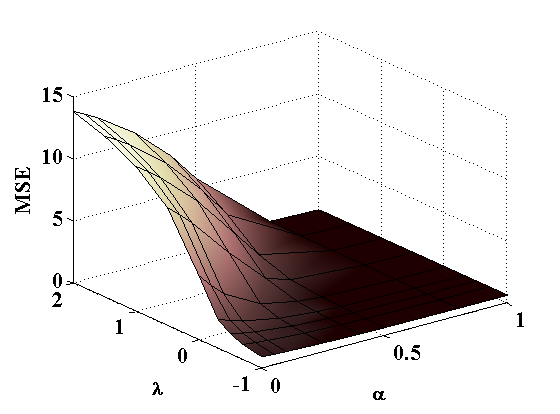}
		\label{FIG:sd_2_20}}
	\caption{Bias and MSE of the  estimator of scale parameter $\sigma$ with different contamination proportion $\epsilon$ for case (ii)}
	\label{FIG:sd_2}
\end{figure}

\begin{figure}[!th]
	\centering
	\includegraphics[width=0.3\textwidth]{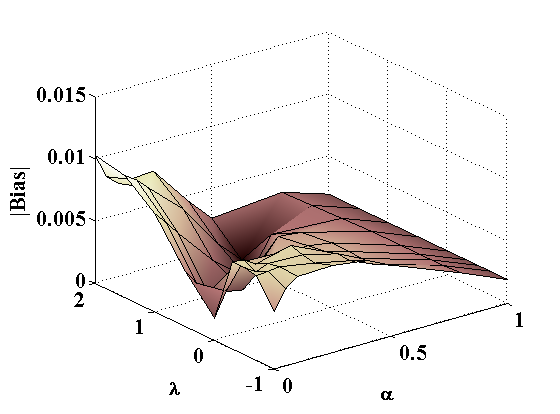}
	\includegraphics[width=0.3\textwidth]{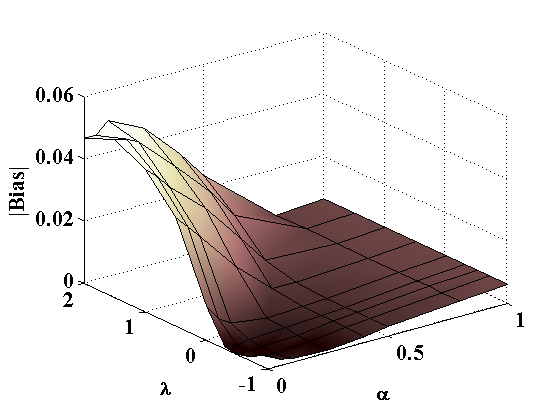}
	\includegraphics[width=0.3\textwidth]{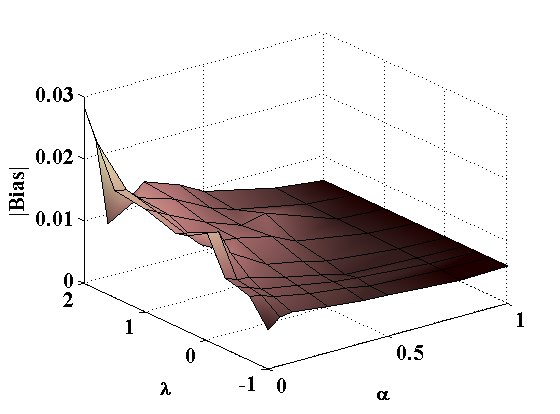}
	\\
	\subfloat[$\epsilon=5\%$]{
		\includegraphics[width=0.3\textwidth]{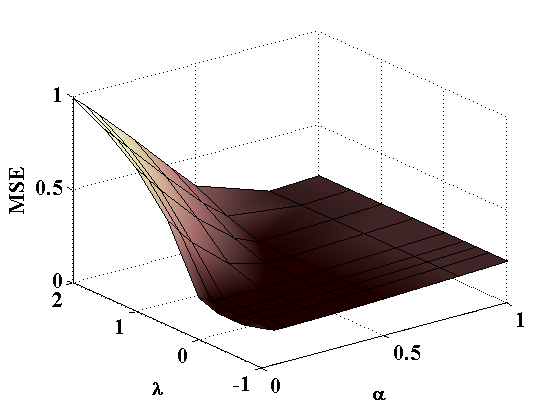}
		\label{FIG:mu_4_5}}
	\subfloat[$\epsilon=10\%$]{
		\includegraphics[width=0.3\textwidth]{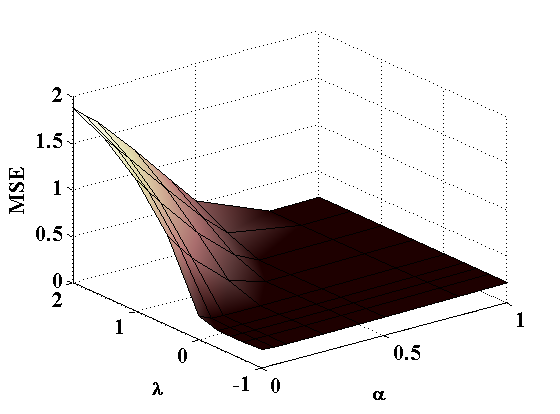}
		\label{FIG:mu_4_10}}
	\subfloat[$\epsilon=20\%$]{
		\includegraphics[width=0.3\textwidth]{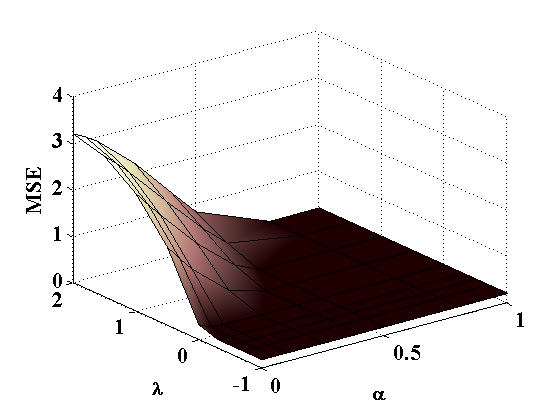}
		\label{FIG:mu_4_20}}
	\caption{Bias and MSE of the  estimator of location parameter $\mu$ with different contamination proportion $\epsilon$ for case (iv)}
	\label{FIG:mu_4}
\end{figure}

\begin{figure}[!th]
	\centering
	\includegraphics[width=0.3\textwidth]{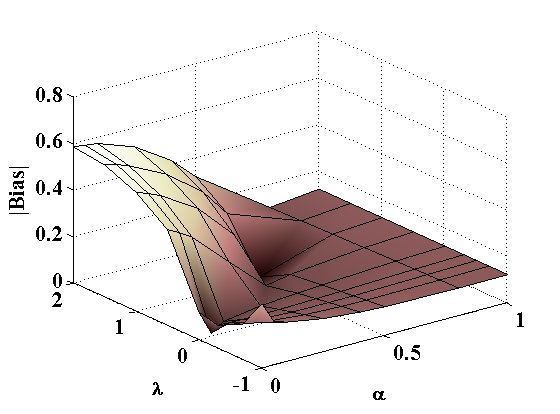}
	\includegraphics[width=0.3\textwidth]{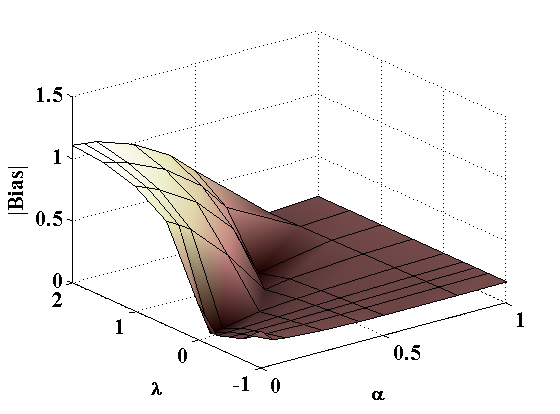}
	\includegraphics[width=0.3\textwidth]{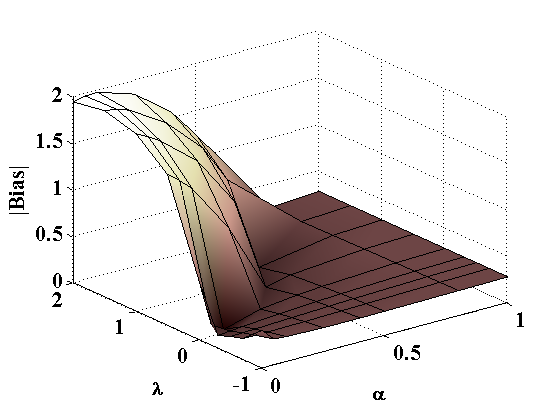}
	\\
	\subfloat[$\epsilon=5\%$]{
		\includegraphics[width=0.3\textwidth]{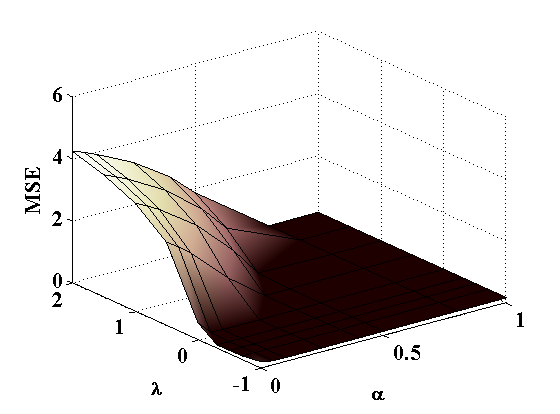}
		\label{FIG:sd_4_5}}
	\subfloat[$\epsilon=10\%$]{
		\includegraphics[width=0.3\textwidth]{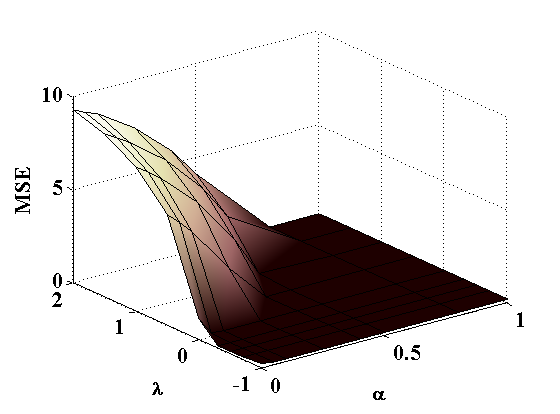}
		\label{FIG:sd_4_10}}
	\subfloat[$\epsilon=20\%$]{
		\includegraphics[width=0.3\textwidth]{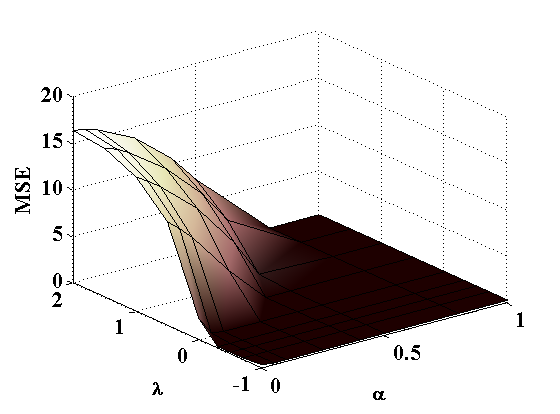}
		\label{FIG:sd_4_20}}
	\caption{Bias and MSE of the  estimator of scale parameter $\sigma$ with different contamination proportion $\epsilon$ for case (iv)}
	\label{FIG:sd_4}
\end{figure}

\begin{figure}[!th]
	\centering
	\includegraphics[width=0.3\textwidth]{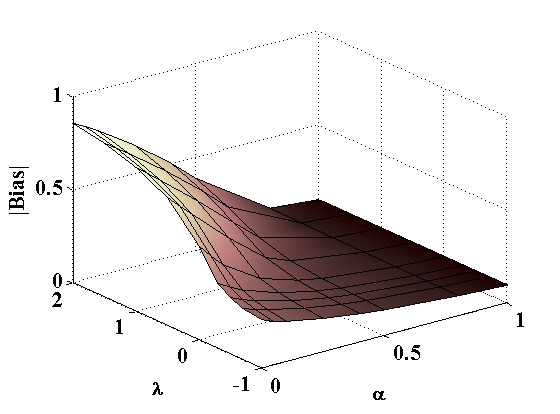}
	\includegraphics[width=0.3\textwidth]{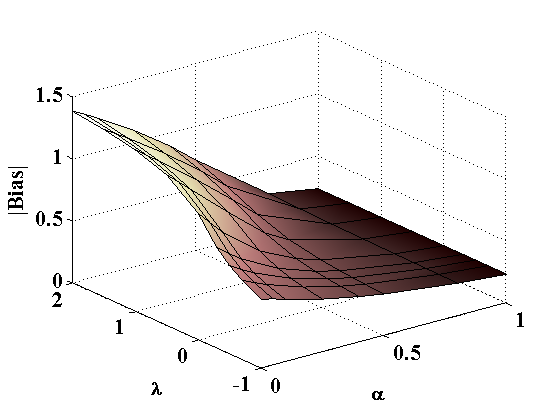}
	\includegraphics[width=0.3\textwidth]{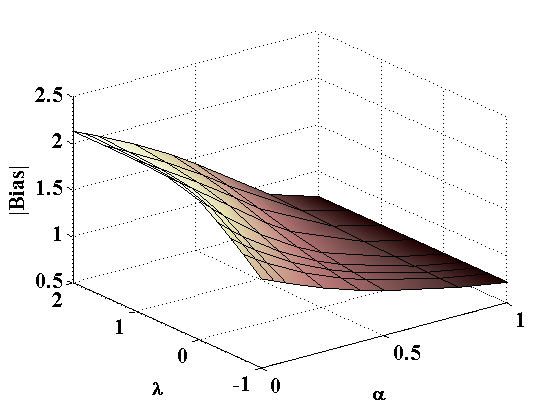}
	\\
	\subfloat[$\epsilon=5\%$]{
		\includegraphics[width=0.3\textwidth]{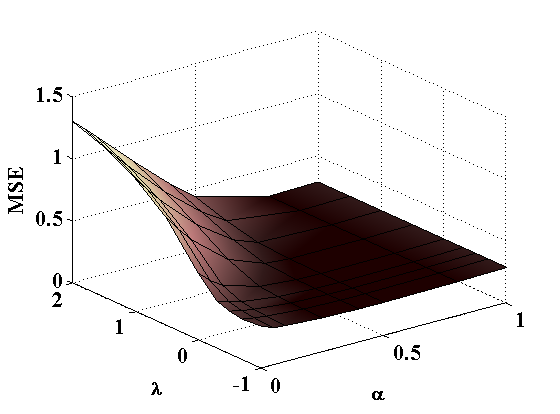}
		\label{FIG:mu_5_5}}
	\subfloat[$\epsilon=10\%$]{
		\includegraphics[width=0.3\textwidth]{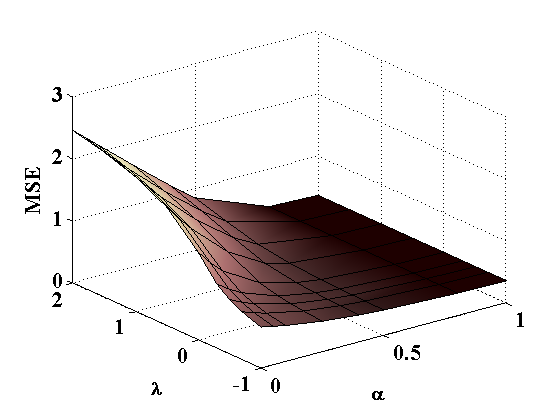}
		\label{FIG:mu_5_10}}
	\subfloat[$\epsilon=20\%$]{
		\includegraphics[width=0.3\textwidth]{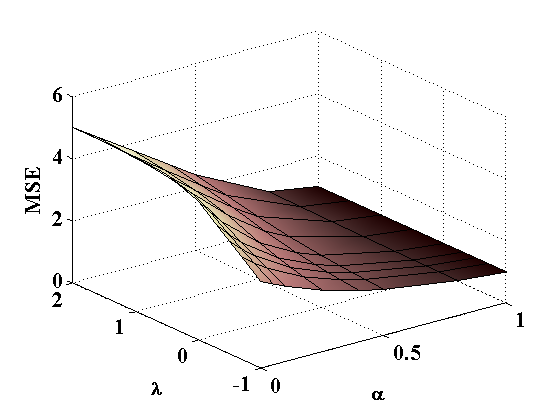}
		\label{FIG:mu_5_20}}
	\caption{Bias and MSE of the  estimator of location parameter $\mu$ with different contamination proportion $\epsilon$ for case (v)}
	\label{FIG:mu_5}
\end{figure}

\begin{figure}[!th]
	\centering
	\includegraphics[width=0.3\textwidth]{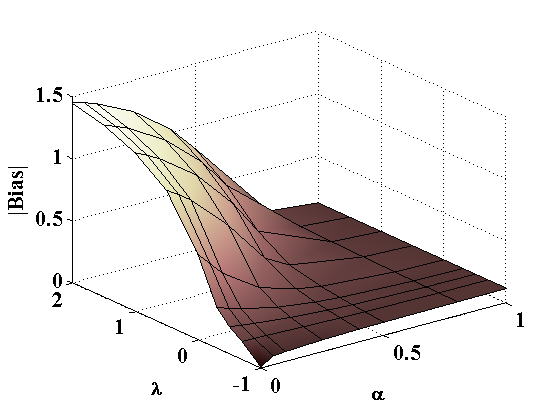}
	\includegraphics[width=0.3\textwidth]{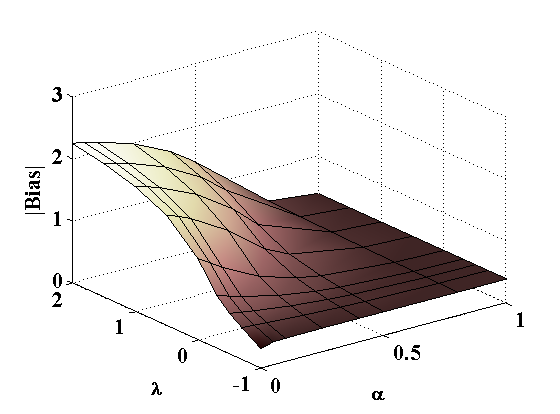}
	\includegraphics[width=0.3\textwidth]{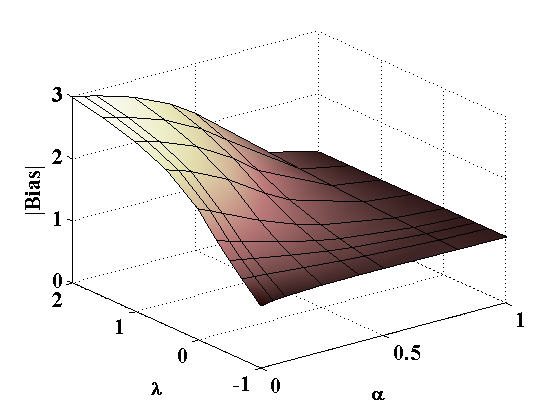}
	\\
	\subfloat[$\epsilon=5\%$]{
		\includegraphics[width=0.3\textwidth]{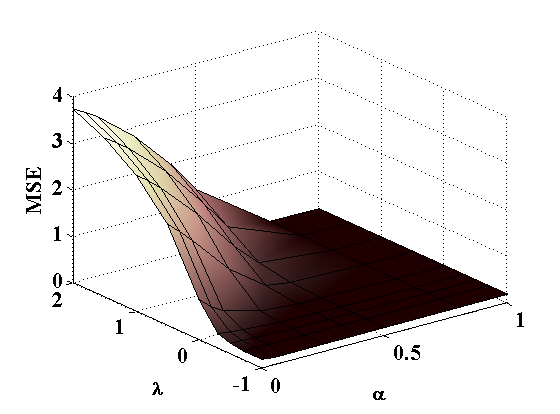}
		\label{FIG:sd_5_5}}
	\subfloat[$\epsilon=10\%$]{
		\includegraphics[width=0.3\textwidth]{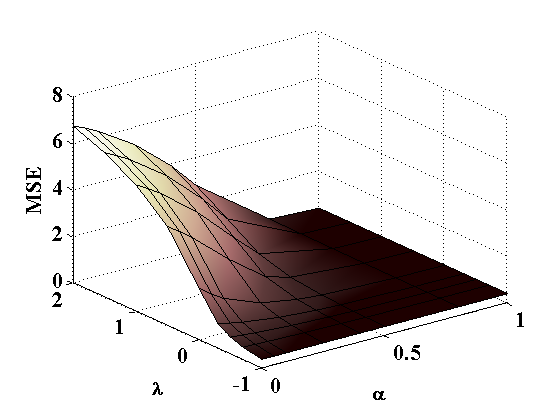}
		\label{FIG:sd_5_10}}
	\subfloat[$\epsilon=20\%$]{
		\includegraphics[width=0.3\textwidth]{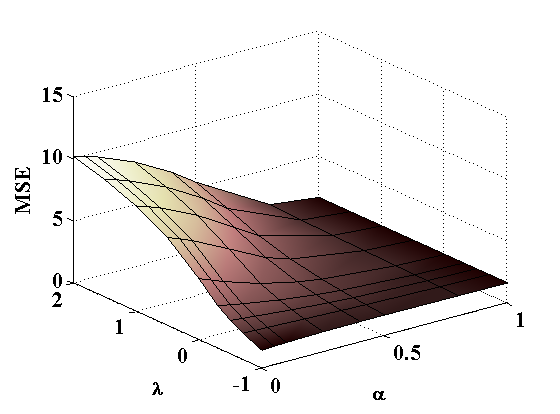}
		\label{FIG:sd_5_20}}
	\caption{Bias and MSE of the  estimator of scale parameter $\sigma$ with different contamination proportion $\epsilon$ for case (v)}
	\label{FIG:sd_5}
\end{figure}

For brevity, we will only present the MSE under some interesting cases of contaminations 
in Figures \ref{FIG:mu_1} to \ref{FIG:sd_5}. If we consider contamination by the
mean shifted density (case (i)), then the variance of the overall distribution also gets affected 
and the MLE of both $\mu$ and $\sigma$ gets affected even under a small contamination of $5\%$. 
The minimum $S^*$-divergence estimator with positive $\lambda$ and small $\alpha$ perform worse than the MLE
under such contaminations; but that corresponding estimators at $\lambda <0$ perform quite robustly ignoring 
the effect of contamination. Even under the heavy $20\%$ contamination in case (i),  
the members of the minimum $S^*$-divergence family with negative $\lambda$ and 
moderately large $\alpha$ successfully ignore the contamination to give 
the robust estimates of the location parameter $\mu$.




%
%

If we consider the contamination by a density with larger variance (as in case (ii)), it should ideally 
affect the scale parameter only -- not the location parameter. The minimum $S^*$-divergence estimators
with negative $\lambda$ and larger $\alpha$ give us quite robust estimates of the scale parameter $\sigma$ 
even with $20\%$ contamination in this case; also the MSE of the location parameter $\mu$ remains almost 
equal to that of the MLE since they are not significantly affected by the contamination in variance. But the minimum 
$S^*$-divergence estimators of $\mu$ with $\lambda>0$ and small $\alpha$ close to zero are seen to be 
affected by such contamination also implying their inferiority compared to the MLE. 
Under the contamination with a symmetric heavy-tail density (case (iv)), 
the performance of the minimum $S^*$-divergence estimators is exactly similar 
to the case of large variance contamination (case (ii));
and their performance under the non-symmetric contamination (case (v)) is again similar to that of 
the mean-shifted contamination (case (i)).
In case (iii) where the true distribution is contaminated by a distribution with smaller variance 
but the same mean, we did not get a consistent pattern to make specific recommendations about 
which choices of the tuning parameter would lead to better estimators in terms of 
their stability under contamination.
As such we have not reported the findings for this case.

Therefore, the overall performance of the minimum $S^*$-divergence estimators can be characterized into 
two groups; just as in the case of minimum $S$-divergence estimators. One group consists of the 
estimators corresponding to $\lambda<0$ and moderate $\alpha$ close to $0.5$ generating highly robust 
estimator with comparable efficiencies with respect to the maximum likelihood estimators. The second group 
consists of those with $\lambda>0$ and small $\alpha$ close to zero which perform even worse than the 
maximum likelihood estimator in terms of the robustness under any type of contamination (except case (iii)). 
The robustness of the minimum $S^*$-divergence estimators are, thus, 
dependent on both the parameters $\lambda$ and $\alpha$.


\section{Real Data Examples}\label{SEC:MS*DE_real data}


\subsection{Short's Data}

\begin{figure}[!fb]
	\centering
	\includegraphics[width=120mm, height=70mm] {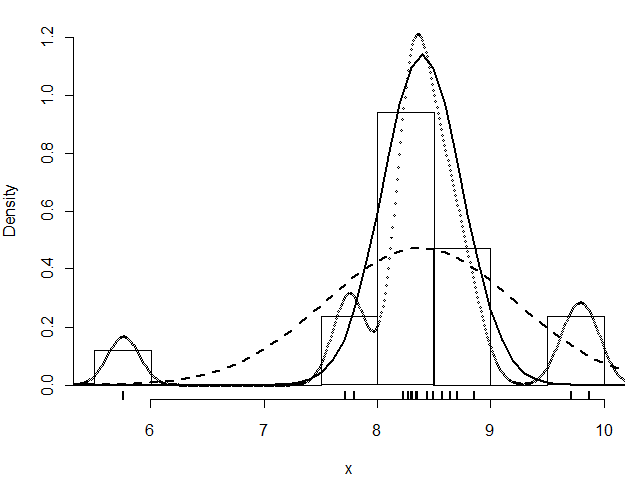}
	\caption{Normal density fits to Short's data (solid line corresponds to $\alpha = 0.5$ and $\lambda=-0.5$, 
		dotted line corresponds to kernel density, dashed line corresponds to MLE)}
\end{figure}

In this example we consider Short's data  for the determination 
of the parallax of the sun, the angle subtended by the earth's radius, as if viewed and measured from 
the surface of the sun. From this angle and available knowledge of the physical dimensions of the earth,
one can easily calculate the mean distance of the earth to the sun. The raw observations are presented in Data 
Set 2 of Stigler (1977).

The data were previously analyzed by many authors including (Basu et al., 2011). 
Observing the pattern of the data, one might model it using a normal density with some 
mean $\mu$ and variance $\sigma^2$. It was observed that 
the data contain a large outlier at 5.76 which severely affects the maximum likelihood estimator.
The maximum likelihood estimates of location $(\mu)$ and scale $(\sigma)$ are 8.378 and 0.846,
respectively. But,  removing the large outlier at 5.76, the maximum likelihood estimates 
of the location and scale become 8.541 and 0.552. So, there is a clear need of using a suitable robust 
technique to estimate the parameters based on this data set.

\begin{table}[!th]
	\centering 
	\caption{Estimates of the location parameter $\mu$ for Short's data. }
	\begin{tabular}{r|lllllll} \hline\noalign{\smallskip}
		$\lambda$	&	$\alpha = $ 0	&	$\alpha = $ 0.1	&	$\alpha = $ 0.3	&	$\alpha = $ 0.4	&	$\alpha = $ 0.5	&	$\alpha = $ 0.8	&	$\alpha = $ 1	\\ \noalign{\smallskip}\hline\noalign{\smallskip}
		$-1$	&	8.38	&	8.39	&	8.39	&	8.39	&	8.40	&	8.40	&	8.41	\\
		$-0.7$	&	8.40	&	8.39	&	8.39	&	8.40	&	8.40	&	8.40	&	8.41	\\
		$-0.5$	&	8.45	&	8.42	&	8.40	&	8.40	&	8.40	&	8.40	&	8.41	\\
		$-0.3$	&	8.49	&	8.47	&	8.42	&	8.41	&	8.40	&	8.40	&	8.41	\\
		0	&	8.38	&	8.43	&	8.46	&	8.43	&	8.41	&	8.40	&	8.41	\\
		0.5	&	8.26	&	8.30	&	8.39	&	8.44	&	8.45	&	8.41	&	8.41	\\
		1	&	8.21	&	8.24	&	8.32	&	8.36	&	8.42	&	8.41	&	8.41	\\
		1.5	&	8.17	&	8.20	&	8.26	&	8.31	&	8.36	&	8.41	&	8.41	\\
		2	&	8.14	&	8.16	&	8.22	&	8.26	&	8.31	&	8.42	&	8.41	\\
		\noalign{\smallskip}\hline
	\end{tabular}
	\label{TAB:mu_short}
\end{table}

\begin{table}[!th]
	\centering
	\caption{Estimates of the scale parameter $\sigma$ for Short's data. }
	\begin{tabular}{r|lllllll} \hline\noalign{\smallskip}
		$\lambda$	&	$\alpha = $ 0	&	$\alpha = $ 0.1	&	$\alpha = $ 0.3	&	$\alpha = $ 0.4	&	$\alpha = $ 0.5	&	$\alpha = $ 0.8	&	$\alpha = $ 1	\\ \noalign{\smallskip}\hline\noalign{\smallskip}
		$-1$	&	0.32	&	0.33	&	0.33	&	0.33	&	0.33	&	0.34	&	0.35	\\
		$-0.7$	&	0.36	&	0.35	&	0.34	&	0.34	&	0.34	&	0.35	&	0.35	\\
		$-0.5$	&	0.45	&	0.41	&	0.36	&	0.35	&	0.35	&	0.35	&	0.35	\\
		$-0.3$	&	0.57	&	0.51	&	0.41	&	0.38	&	0.36	&	0.35	&	0.35	\\
		0	&	0.85	&	0.76	&	0.51	&	0.44	&	0.39	&	0.35	&	0.35	\\
		0.5	&	0.97	&	0.93	&	0.79	&	0.65	&	0.49	&	0.36	&	0.35	\\
		1	&	1.01	&	0.98	&	0.89	&	0.81	&	0.69	&	0.37	&	0.35	\\
		1.5	&	1.03	&	1.01	&	0.94	&	0.88	&	0.80	&	0.38	&	0.35	\\
		2	&	1.04	&	1.02	&	0.97	&	0.92	&	0.86	&	0.39	&	0.35	\\
		\noalign{\smallskip}\hline
	\end{tabular}
	\label{TAB:sigma_short}
\end{table}

We will apply the proposed minimum $S^*$-divergence estimators here with a normal kernel $W(x, y, h)$ 
as discussed in the previous sections. The bandwidth will be chosen from the normal reference rule 
as in the case of  our simulation studies, described in Section \ref{SEC:4MS*DE_simulation}.
In Tables \ref{TAB:mu_short} and \ref{TAB:sigma_short}, we provide the minimum $S^*$-divergence 
estimates of the location and scale parameters under several values of $\alpha$ and $\lambda$. 
We can clearly see that many of the $S^*$-divergence  estimators are not 
affected by the outliers. This includes the $S^*$-divergence with negative $\lambda$ with small
$\alpha$ and the same with positive $\lambda$ but large positive $\alpha$.
However the $S^*$-divergence  estimators with $\lambda > 0$ and small positive  $\alpha$ are very 
sensitive with respect to the outlier like the MLE (which corresponds to $\alpha =0$ and $\lambda=0$). 
Thus the performance of the minimum $S^*$-divergence  estimators for the continuous model is also 
exactly similar to that of the minimum $S$-divergence estimators under the discrete model.
There is a clear triangular region, in the lower left hand corner of Table \ref{TAB:sigma_short}, 
where the estimators are adversely affected by the outlier, in some cases severely.  
This is consistent with our findings of Section \ref{SEC:4MS*DE_simulation}.

\subsection{Newcomb's Data}

This example involves Newcomb's light speed data 
(Stigler, 1977, Table 5). The data set shows a nice unimodal structure, and the normal model would
have provided an excellent fit to the data if the two large outliers were not there.
The data set was previously analyzed by many researchers who had shown the need of 
suitable robust estimation methods for fitting a normal model because the usual 
maximum likelihood estimator is seen to be highly affected by the presence of the two outliers in the data set. 
So, we will compute  the robust minimum $S^*$-divergence estimators for the present data, 
with the kernel and the corresponding bandwidth being the same as in the previous example. 
The estimates obtained for different members of the $S^*$-divergence family are presented in 
Tables \ref{TAB:mu_newcomb} and \ref{TAB:sigma_newcomb}.

\begin{table}[ht]
	\centering 
	\caption{Estimates of the location parameter $\mu$ for Newcomb's data. }
	\begin{tabular}{r| lllllll} \hline\noalign{\smallskip}
		$\lambda$	&	$\alpha = $ 0	&	$\alpha = $ 0.1	&	$\alpha = $ 0.3	&	$\alpha = $ 0.4	&	$\alpha = $ 0.5	&	$\alpha = $ 0.8	&	$\alpha = $ 1	\\ \noalign{\smallskip} \hline \noalign{\smallskip}
		$-1$	&	27.72	&	27.69	&	27.63	&	27.60	&	27.56	&	27.47	&	27.42	\\
		$-0.7$	&	27.73	&	27.70	&	27.63	&	27.60	&	27.57	&	27.47	&	27.42	\\
		$-0.5$	&	27.73	&	27.70	&	27.64	&	27.60	&	27.57	&	27.47	&	27.42	\\
		$-0.3$	&	27.72	&	27.70	&	27.64	&	27.61	&	27.57	&	27.47	&	27.42	\\
		0	&	26.21	&	27.58	&	27.64	&	27.61	&	27.57	&	27.47	&	27.42	\\
		0.5	&	22.69	&	23.80	&	26.47	&	27.57	&	27.57	&	27.48	&	27.42	\\
		1	&	21.39	&	22.23	&	24.25	&	25.52	&	27.27	&	27.48	&	27.42	\\
		1.5	&	20.62	&	21.29	&	22.97	&	24.06	&	25.41	&	27.48	&	27.42	\\
		2	&	20.12	&	20.66	&	22.09	&	23.05	&	24.25	&	27.48	&	27.42	\\
		\noalign{\smallskip}\hline
	\end{tabular}
	\label{TAB:mu_newcomb}
\end{table}

\begin{table}[h]
	\centering
	\caption{Estimates of the scale parameter $\sigma$ for Newcomb's data. }
	\begin{tabular}{r| lllllll} \hline\noalign{\smallskip}
		$\lambda$	&	$\alpha = $ 0	&	$\alpha = $ 0.1	&	$\alpha = $ 0.3	&	$\alpha = $ 0.4	&	$\alpha = $ 0.5	&	$\alpha = $ 0.8	&	$\alpha = $ 1	\\ \noalign{\smallskip}\hline\noalign{\smallskip}
		$-1$	&	4.94	&	4.97	&	4.98	&	4.97	&	4.95	&	4.90	&	4.87	\\
		$-0.7$	&	4.99	&	5.00	&	4.99	&	4.98	&	4.96	&	4.90	&	4.87	\\
		$-0.5$	&	5.02	&	5.02	&	5.00	&	4.98	&	4.96	&	4.90	&	4.87	\\
		$-0.3$	&	5.10	&	5.06	&	5.01	&	4.99	&	4.97	&	4.90	&	4.87	\\
		0	&	10.66	&	5.48	&	5.05	&	5.01	&	4.98	&	4.90	&	4.87	\\
		0.5	&	18.63	&	17.10	&	10.76	&	5.19	&	5.01	&	4.91	&	4.87	\\
		1	&	20.48	&	19.53	&	16.49	&	13.69	&	7.02	&	4.91	&	4.87	\\
		1.5	&	21.36	&	20.65	&	18.51	&	16.75	&	13.89	&	4.92	&	4.87	\\
		2	&	21.88	&	21.30	&	19.61	&	18.29	&	16.29	&	4.92	&	4.87	\\
	\noalign{\smallskip}	\hline
	\end{tabular}
	\label{TAB:sigma_newcomb}
\end{table}

Many of the minimum $S^*$-divergence estimates automatically discount these large observations, 
unlike the maximum likelihood estimate (corresponding to $\alpha=0$ and $\lambda=0$ in 
the $S^*$-divergence measure).  
Once again the minimum $S^*$-divergence estimates for all $\alpha$ with $\lambda<0$ and for large $\alpha$ 
close to one with $\lambda \ge 0$  are remarkably close to each other and 
robust with respect to the large outliers.
Once again the region of instability is in the lower left hand corner of the tables.

\section{A Suggestion for the Kernels : $\alpha$-Transparent Kernel}
\label{SEC:4Trans_kernel}

In this paper, we have developed the properties of the minimum $S^*$-divergence estimators 
under the continuous models following the approach of Basu and Lindsay (1994).
The work of Basu and Lindsay (1994) was in the context of minimum disparity estimators,
where they have also derived some conditions on the kernel density so that 
the minimum disparity estimators obtained by the smoothed model approach under continuous models 
have similar asymptotic distribution as that of the original minimum disparity estimators 
without any kernel smoothing under discrete models and hence they are fully efficient.
They have termed those special kernels as the ``Transparent kernel",
although there are only few examples of them. In this section, we will derive 
similar conditions on the kernel so that the asymptotic distribution of the proposed 
minimum $S^*$-divergence estimators under continuous models (as derived in previous sections) 
have the  same form as that of the minimum $S$-divergence estimators under discrete models  
(as derived in Ghosh, 2014).

Let us again start with the special case corresponding to $\lambda = 0$ 
which gives the well-known density power divergence; in this case we can compare 
the results of the corresponding smoothed model based estimators
with the existing results of the minimum density power divergence estimators.
In particular, comparing the estimating equation (\ref{EQ:est-eqn-MDPD_star}) with 
the MDPDE estimating equation (\ref{EQ:est-eqn-MDPD}) we can find the required conditions on the 
kernel function under which the $\textrm{MDPDE}^*$ (the $\textrm{MSDE}^*$ for $\lambda = 0$) 
coincides with the MDPDE. 
The result is presented in the following Lemma.

\bigskip
\begin{lemma}\label{LEM:4MDPD*_trans}
	Suppose the kernel function $W(x,y,h)$ used in the smoothing densities is such that 
	\begin{eqnarray}\label{EQ:Trans_Kernel_01}
	u_\theta^{\alpha*}(y) = M f_\theta^\alpha(y) {u_\theta}(y) + L
	\end{eqnarray} 
	for a $p$-vector $L$ depending only on $\alpha$ and $h$,
	and a $p\times p$ nonsingular matrix $M$ depending on $\theta$, $\alpha$ and $h$  
	where for each components $\theta_j$ of $theta$ ($j=1,\ldots,p$), we have
	\begin{equation}
	\mbox{either ``$\int u_{\theta_j} f_{\theta}^{1+\alpha} = 0"$, or 
		``the $j$-th column of $M$ is independent of $\theta$".}
	\label{EQ:Trans_Kernel_02}
	\end{equation}  
	Then the estimating equation 
	for the Minimum DPD$^*$ Estimator is the same as that for the Minimum DPD Estimator and hence
	the two estimators are indeed equal.\hfill{$\square$}
\end{lemma}

Further, just as in the case of disparities, 
the condition (\ref{EQ:Trans_Kernel_01}) imposed on the kernel function 
also ensure that the asymptotic variance of the $\textrm{MSDE}^*$s will be equal to that
of the MSDEs under discrete model, beyond this special case of $\lambda =0$ also.
We will justify this condition further by comparing the influence functions and asymptotic distributions of 
the two types of estimators in the following two corollaries.

\bigskip
\begin{corollary}\label{COR:4IF_MSDE*_model_Trans}
	Suppose the true density $g$ belongs to the model family  
	$ \{ f_\theta : \theta \in \Theta \}$,  i.e., $ g = f_\theta $ for some $\theta\in\Theta$.
	Also assume that the kernel function $W(x,y,h)$ used in the smoothing 
	satisfies the condition (\ref{EQ:Trans_Kernel_01}) for some nonsingular matrix $M$ and vector $L$
	as in Lemma \ref{LEM:4MDPD*_trans} along with Condition (\ref{EQ:Trans_Kernel_02}). 
	Then the Influence Function for the  minimum $S^*$-divergence estimator functional 
	$T_{(\alpha, \lambda)}^*$ at the distribution $G = F_\theta$  becomes 
	the same as that of the minimum $S$-divergence estimator at $G = F_\theta$, 
	given in Equation (\ref{EQ:S_div_IF_model}).
\end{corollary}
\noindent\textbf{Proof: }
In view of  Corollary \ref{COR:4IF_MSDE*_model} and Equation (\ref{EQ:4cont_J*_model2}),
we only need to show that, Condition (\ref{EQ:Trans_Kernel_01}) 
along with (\ref{EQ:Trans_Kernel_02}) implies 
\begin{eqnarray}
E_\theta [ - \nabla u_{\theta}^{\alpha*}(X) ] = M \int u_\theta(x)u_\theta(x)^Tf_\theta^{1+\alpha}(x)dx.
\label{EQ:4cont_nabla_u^*}
\end{eqnarray}
Then, it follows that 
\begin{eqnarray}
IF(y;F_\theta,T_{(\alpha, \lambda)}^*) &=& [J^{*}(\theta)]^{-1} \left\{ u_{\theta}^{\alpha*}(y) - 
E_\theta [ u_{\theta}^{\alpha*}(X) ] \right\}\nonumber\\
&=& \left[M \int u_\theta(x)u_\theta(x)^Tf_\theta^{1+\alpha}(x)dx\right]^{-1} \nonumber\\
&& ~~~~ \times \left\{M u_\theta(y)f_\theta^\alpha(y) 
+ L - E_\theta [ M u_\theta(X)f_\theta^\alpha(X) + L] \right\} \nonumber\\
&=& \left[\int u_\theta(x)u_\theta(x)^Tf_\theta^{1+\alpha}(x)dx\right]^{-1} \nonumber\\
&& ~~~~\times\left\{ u_\theta(y)f_\theta^\alpha(y) - E_\theta [ u_\theta(X)f_\theta^\alpha(X)] \right\}, 
\nonumber
\end{eqnarray}
which is the same as given in (\ref{EQ:S_div_IF_model}); it proves the corollary.

Now, to prove (\ref{EQ:4cont_nabla_u^*}), let us take the derivative with respect to $\theta$ 
on both sides of the condition (\ref{EQ:Trans_Kernel_01}) to get
\begin{eqnarray}
\nabla u_{\theta}^{\alpha*}(x)&=& M \nabla u_\theta (x)f_\theta^{\alpha}(x) 
+ \alpha M u_\theta(x)u_\theta(x)^Tf_\theta^{\alpha}(x) \nonumber\\
&& + \left[(\nabla_1M)u_\theta(x) ~~ (\nabla_2M)u_\theta(x) ~~ \cdots~~ 
(\nabla_pM)u_\theta(x)\right]f_\theta^{\alpha}(x), ~~~~
\label{EQ:4cont_eq0}
\end{eqnarray}
where $\nabla_j$ represents the derivative with respect to the $j$-th component of $\theta$.
Taking expectation with respect to $f_\theta$ in both sides of (\ref{EQ:4cont_eq0}), we have
\begin{eqnarray}
E_\theta [ - \nabla u_{\theta}^{\alpha*}(X) ] &=& M \int \nabla u_\theta(x)f_\theta^{1+\alpha}(x)dx \nonumber\\
&& ~~~~~~~~~+ \alpha M \int u_\theta(x)u_\theta(x)^Tf_\theta^{1+\alpha}(x)dx,~~~~~~
\label{EQ:4cont_eq1}
\end{eqnarray}
since the expectation of the third term in (\ref{EQ:4cont_eq0}) is zero by 
Condition (\ref{EQ:Trans_Kernel_02}).
Also, integrating the first integral in (\ref{EQ:4cont_eq1}) by parts, we get 
\begin{eqnarray}
\int \nabla u_\theta(x)f_\theta^{1+\alpha}(x)dx = 
- (1+\alpha) \int u_\theta(x)u_\theta(x)^Tf_\theta^{1+\alpha}(x)dx.
\label{EQ:4cont_eq2}
\end{eqnarray}
Combining (\ref{EQ:4cont_eq1}) and (\ref{EQ:4cont_eq2}) we get the desired result (\ref{EQ:4cont_nabla_u^*}),
completing the proof.
\hfill{$\square$}

\begin{corollary}\label{COR:4cont_asymp_model2}
	Suppose the true density $g$ belongs to the model family  $\{ f_\theta : \theta \in \Theta \}$,
	i.e., $ g = f_\theta $ for some $\theta \in \Theta$. 
	Also assume that the kernel function $W(x,y,h)$ used in the smoothing 
	satisfies the condition (\ref{EQ:Trans_Kernel_01}) for some nonsingular matrix $M$ and vector $L$
	as in Lemma \ref{LEM:4MDPD*_trans} along with Condition (\ref{EQ:Trans_Kernel_02}). 
	Then the asymptotic distribution for the  minimum $S^*$-divergence estimator 
	$\theta_n^*$ is normal with mean zero and variance-covariance matrix as given by 
	the expression (\ref{EQ:4cont_asymp_var_model}).
	\hfill{$\square$}
\end{corollary}

Thus we have seen that if the kernel function satisfies Equation (\ref{EQ:Trans_Kernel_01})  
then the  minimum $S^*$-divergence estimator and the  minimum $S$-divergence estimator have 
the same influence function at the model distribution.
Thus, under this assumption, all the first order asymptotic properties of the $\textrm{MSDE}^*$s 
should be same as that of the MSDEs under discrete models. 
In particular, we have also seen such equivalence in terms of their asymptotic distribution.
We will refer to kernel function satisfying  Equation (\ref{EQ:Trans_Kernel_01}) 
as the \textit{$\alpha$-transparent kernel} for the model $f_\theta$;
at $\alpha=0$ this notion coincides with that of the ``transparent kernel" 
as defined in Basu and Lindsay (1994).
\index{$\alpha$-transparent Kernel}

\begin{definition}\label{DFN:alpha_transparent_kernel}
	Consider the parametric model ${\cal F} = \{F_\theta: \theta \in \Theta \subseteq \mathbb{R}^p\}$. 
	Let $u_\theta^{\alpha*}(x)$ be as defined in Equation (\ref{EQ:u_theta_alpha_star}).
	Then a kernel function $W(x, y, h)$ will be called a $\alpha$-transparent kernel for the above family
	of models if it satisfies the Condition (\ref{EQ:Trans_Kernel_01})
	for some nonsingular matrix $M$ and vector $L$ as in Lemma \ref{LEM:4MDPD*_trans}
	along with Condition (\ref{EQ:Trans_Kernel_02}).
\end{definition}

For the mean parameter of the normal model, the Gaussian kernel provides 
one example of an $\alpha$-transparent kernel at $\alpha = 0$; 
in this case the asymptotic variance of the MSDE$^*$ becomes independent of the bandwidth $h$ 
as seen in Remark \ref{REM:4asymp_var_H0}. Although the calculations are not provided here, 
the same is true for the normal variance. In general, however, the $\alpha$-transparent kernel 
is a theoretical construct and we do not have other examples of $\alpha$-transparent kernels at this point. 
If one does exist, the problem under consideration will have simple solutions. 
However, even if one does not exist, we expect that letting $h\rightarrow 0$ will allow 
the asymptotic variance of the MSDE$^*$ to stabilize around the expression in 
equation (\ref{EQ:4cont_asymp_var_model}), as discussed in Section \ref{SEC:MS*DE_asymptotic}.

\section{Role of $\lambda$ in Robustness: Second Order Influence Analysis}\label{SEC:Second_order}

We have seen, in all our empirical studies, that the robustness of 
the minimum $S^*$-divergence estimator depends crucially on
both the tuning parameters, although no theoretical evidence is discovered as 
to the role of $\lambda$ through the classical (first order)
influence function analysis presented in Section \ref{SEC:MS*DE_IF}. 
Similar limitations of the classical influence function 
in describing the robustness of the minimum divergence estimators have also been observed by Lindsay (1994) 
and Ghosh et al.~(2013) for the disparity family and the $S$-divergence family under the discrete set-up; 
these authors have suggested that a second order influence function analysis of the minimum divergence estimator 
can provide us a better indication of the theoretical robustness.  
In this Section, we briefly present the second order influence function analysis 
for the minimum $S^*$-divergence estimator to examine the role of $\lambda$ in the context of robustness.

The influence function of any statistical functional provides 
a first order approximation to the bias caused by contamination.
In case of the MSDE$^*$, the first order Taylor series approximation yields  
$$
\Delta T_{(\alpha,\lambda)}^*(\epsilon) = T_{(\alpha,\lambda)}^*(G_\epsilon ) - T_{(\alpha,\lambda)}^*( G ) 
\approx \epsilon T'(y),
$$
where $T'(y) = IF(y; G, T_{(\alpha, \lambda)}^*)$ and $G_\epsilon$ is as defined in Section \ref{SEC:MS*DE_IF}. 
Therefore the predicted bias up to the first order will be the same for all $\lambda$,
which provides and insufficient description of the robustness of the MSDE$^*$s. 
So, we consider the second order Taylor series expansion 
of the predicted bias to get 
$$\Delta T_{(\alpha,\lambda)}(\epsilon) = \epsilon T'(y) + \frac{\epsilon^2}{2} T''(y),$$
where 
$T''(y) = \frac{\partial}{\partial\epsilon}T_{(\alpha,\lambda)}(G_\epsilon)\big|_{\epsilon=0}$.
See Ghosh et al.~(2013) for a more detailed discussion.
In the following theorem, we present the expression of this second order approximation 
$T''(y)$ for the case of scalar parameter; the proof is elementary and hence omitted. 
Multi-parameter extensions can be done in a straightforward manner.

\begin{theorem}
Consider the set-up of Section \ref{SEC:MS*DE_IF} with $\theta$ being a scalar parameter. 
Assume that the true distribution belongs to the model $\mathcal{F}$. 
For the minimum $S^*$-divergence estimator, we have 
$$
T''(y) = T'(y) [J^*]^{-1} \left[m_1^*(y) + \lambda(1-\alpha)m_2^*(y)\right],
$$ 
where  
\begin{eqnarray}
m_1^*(y)& =&  2u_\theta^{1\alpha*}(y) + 2\alpha u_\theta^{2\alpha\ast}(y) - 2J^* 
- T'(y) \left[ (1+2\alpha)D_1^* + 3D_2^*  \right], \nonumber  \\
m_2^*(y) &=& D_1^* T'(y) - 2 u_\theta^{2\alpha\ast}(y) 
+ \frac{u_\theta^{(\alpha-1)\ast}(y) - E_\theta [ u_{\theta}^{\alpha*}(X) ]}{
u_{\theta}^{\alpha*}(y) -E_\theta [ u_{\theta}^{\alpha*}(X) ]}. \nonumber  
\end{eqnarray}
with $D_1^*=\int \widetilde{u_\theta}^3 (f_\theta^*)^{1+\alpha}$,
$D_2^*=\int \widetilde{u_\theta} (\nabla \widetilde{u_\theta}) (f_\theta^*)^{1+\alpha}$ and 
$$u_\theta^{1\alpha*}(y) = \int \nabla \widetilde{u_\theta}(x)\{f_\theta^*(x)\}^\alpha W(x,y,h)dx. $$
\end{theorem}

Note that the second order influence function of the MSDE$^*$ depends on the tuning parameter $\lambda$
through the function $m_2^*(y)$. The effect of $\lambda$ diminishes as $\alpha \rightarrow 1$
as is expected from the special structure of the $S$-divergence measure; the measure becomes independent of 
$\lambda$ at $\alpha=1$. Further it can be seen that, for most parametric models,
the first order influence function $T'(y)$ and the function $m_1^*$ is bounded for all $\alpha>0$
and the function $m_2^*(y)$ contains 
a term of the order of $\frac{1}{f_\theta^*(y)}$ which can potentially be unbounded. 
Thus whenever $y$ is far from the data cloud one should keep the factor $\lambda(1-\alpha)$ 
small so that the effect of the terms $m_2^*(y)$ is minimized;
this mandates that $\lambda$ should be close to zero and $\alpha$ should be close to 1.
In fact small negative values of $\lambda$ may be preferred as it help to counter balance
the effect of the first term.
As an illustration, we will consider the case of the normal mean parameter.

\begin{remark}
Consider the normal model $N(\theta, \sigma^2)$ with unknown mean $\theta$ and known variance $\sigma^2$. 
Then simple calculations yield $J^* = (2\pi)^{\alpha/2}(1+\alpha)^{3/2}(\sigma^2+h^2)^{-\frac{\alpha+2}{2}}$,
$D_1^*=D_2^*=0=E_\theta [u_{\theta}^{\alpha*}(X)]$ and
\begin{eqnarray}
u_\theta^{\alpha\ast}(y) &=& C_{\alpha,h}(y-\theta)e^{-\frac{\alpha(y-\theta)^2}{2(\alpha h^2+h^2+\sigma^2)}},
\nonumber\\
u_\theta^{2\alpha\ast}(y) &=& C_{\alpha,h}\left\{\frac{h^2}{h^2+\sigma^2}
+\frac{(y-\theta)^2}{(\alpha h^2+h^2+\sigma^2)}\right\}
e^{-\frac{\alpha(y-\theta)^2}{2(\alpha h^2+h^2+\sigma^2)}},\nonumber\\
u_\theta^{1\alpha\ast}(y) &=& C_{\alpha,h}\left(\frac{\alpha h^2}{\sigma^2 + h^2} + 1 \right) 
e^{-\frac{\alpha(y-\theta)^2}{2(\alpha h^2+h^2+\sigma^2)}},
\nonumber
\end{eqnarray}
with $C_{\alpha,h} = (2\pi)^{-\alpha/2}(\sigma^2+h^2)^{-\frac{\alpha-1}{2}}(\alpha h^2+h^2+\sigma^2)^{-3/2}$.
So, we can easily compute the first and second order influence functions of the MSDE$^*$ of $\theta$.
It is easy to see that the first order influence function is bounded for all $\alpha>0$ (independently of $\lambda$)
implying the robustness of the estimator. However, the second order influence function depends on $\lambda$ and
contains a term of the form 
$$\lambda(1-\alpha)Ke^{\frac{(y-\theta)^2}{2(\alpha h^2+h^2+\sigma^2)}},$$
for a constant $K$ depending only on $\alpha$, $h$ and $\sigma^2$. 
So, the second order prediction of robustness of the MSDE$^*$ depends on the parameter $\lambda$
unless $\alpha=1$; the effect of this term decreases as $|\lambda|$ decreases
and  the actual predicted bias (in second order approximation) decreases at negative $\lambda$
counter balancing the effect of the other (bounded) term $m_1^*(y)$.
This fact is quite in-line with our simulation results where a negative $\lambda$ value near zero yields 
the most robust estimator at positive $\alpha$. 
\end{remark}

\section{Choice of the Tuning parameters $\alpha$ and $\lambda$}\label{SEC:Choice_tuning_parameter}

Finally we note that, combining all these asymptotic properties and empirical findings presented above, 
there exists a clear trade-off between efficiency and robustness between different members of the 
$S$-divergence family. This trade-off solely depends on the tuning parameters $\alpha$ and $\lambda$
defining the divergence measure and is exactly similar to that observed for the minimum 
$S$-divergence estimators in discrete cases (Ghosh et al., 2013; Ghosh, 2014).
In particular, efficiency decreases and robustness increases as $\alpha$ increases.
However the first order efficiency and first order influence function are 
theoretically independent of the other tuning parameter $\lambda$, though numerical illustrations 
suggest otherwise -- the estimators with $\lambda<0$, 
and the estimators with $\lambda>0$ with larger values of $\alpha$ 
are more robust than the other members of the family; 
this is consistent with the higher order influence analysis of the previous section.
Therefore, in case of continuous model also 
we can choose an optimum range of values of the tuning parameters following the logic discussed in Ghosh (2014) 
for the discrete case. Combining the asymptotic and empirical findings, it was suggested that
low positive values of $\alpha \in [0.1, 0.25]$ coupled with moderately small negative values of
$\lambda \in [-0.3, -0.5]$ be considered as appropriate choices of the tuning parameter to 
obtain highly robust estimators with minimal loss of efficiency.
We also stick with these choice of tuning parameters for the continuous models and 
the minimum $S^*$-divergence estimators considered here as empirical suggestions.

However, even within this general window there may be a substantial variation in the performance of
the estimators over the tuning parameters, and being able to further pinpoint 
the most suitable estimator in a given practical situation is of importance. 
The applicability of the proposed method will surely be enhanced 
if we can provide data driven ``optimal" choices of the tuning parameter for any given situation.
These kind of approaches are already available in the literature;
most of them are based on some numerical technique like bootstrap, cross-validation or BIC 
(Altman and Leger, 1994; Fang et al., 2013, Park et al., 2014, Kawano, 2014) or 
the minimization of the estimated mean square error (MSE) of the estimator 
(Hong and Kim, 2001; Warwick and Jones, 2005; Ghosh and Basu, 2013b) etc.
We can adopt any such approach in the present context of the minimum $S^*$ divergence estimator to 
select optimum tuning parameters. However, considering the length of the current paper, we only briefly 
mention one such approach, namely the approach taken by Warwick and Jone (2005) with intuitive justification; 
we hope to take up its empirical performance with respect to other competitive approaches  
in a separate paper in the future.

Suppose $\theta^0$ denotes the target parameter and consider a contaminated version of the target density 
$g(x) = (1-\epsilon) f_{\theta^0}(x) + \epsilon \delta_y(x)$ 
in terms of the Dirac delta function $\delta_y(x)$ at the point $y$. 
For given $\alpha$ and $\lambda$, let 
$\theta_{\alpha,\lambda}=\arg\min_{\theta} S_{(\alpha,\lambda)}(g^*,f_\theta^*)$. 
For any given sample dataset from $g$, let $\hat{\theta}_{\alpha,\lambda}$ denote the MSDE$^*$ of $\theta$. 
In the Warwick and  Jones (2005) approach, one minimizes an asymptotic approximation of the summed MSE 
$E\left[ (\hat{\theta}_{\alpha,\lambda} - \theta^0)^T(\hat{\theta}_{\alpha,\lambda} - \theta^0)\right]$ 
to choose the optimal $\alpha$ and $\lambda$.
Under standard regularity conditions as discussed in Section \ref{SEC:MS*DE_asymptotic}, 
we have the approximation 
\begin{equation}
E\left[ (\hat{\theta}_{\alpha,\lambda} - \theta^0)^T(\hat{\theta}_{\alpha,\lambda} - \theta^0)\right] 
=  (\theta_{\alpha,\lambda} - \theta^0)^T(\theta_{\alpha,\lambda} - \theta^0) + 
\frac{1}{n}{\rm trace}\left\{[J^*]^{-1} V^*[J^*]^{-1}\right\},\nonumber
\end{equation}
where $J^*$ and $V^*$ are as in Corollary \ref{COR:4cont_asymp_model1}.
We can easily estimate the asymptotic variance component by substituting  $\hat{\theta}_{\alpha,\lambda}$ 
for ${\theta}_{\alpha,\lambda}$  and $g_n^*$ in place of $g^*$; 
also ${\theta}_{\alpha,\lambda}$ may be estimated by the MSDE$^*$ $\hat{\theta}_{\alpha,\lambda}$. 
But there is no obvious choice of $\theta^0$. As suggested by Warwick and  Jones (2005), 
we may consider several ``pilot" estimators of $\theta$ 
in place of $\theta^0$ in the minimization process and compared its effect through simulation studies. 
Through several simulation in the context of minimum density power divergence estimators (MDPDE), 
they recommended the use of the MDPDE corresponding to $\alpha=1$ as the pilot estimate.

\section{Concluding Remarks }\label{SEC:discussion}

In this work we have developed the theoretical properties of the minimum divergence estimators 
under the general framework of the continuous models. The $S$-divergence family
provides a large collection of divergences with different properties; 
it includes several divergences where the corresponding estimators have near optimal efficiency properties 
and strong robustness properties, and provides a general framework 
to explore the properties of the minimum divergence estimators 
through two tuning parameter $\alpha\geq 0$ and $\lambda\in\mathbb{R}$.
Ghosh et al.~(2013) explored the robustness performances of this general class of minimum divergence estimators 
through the influence function and breakdown point analysis;
Ghosh (2014) proved the asymptotic properties of these estimators under the discrete models
and linked the theoretical properties with the empirical findings and robustness 
presented in Ghosh et al.~(2013). Several other applications of this family of $S$-divergence 
measures are under the lens of many recent researches are also being considered elsewhere; 
see for example, Ghosh, Maji and Basu (2013) and Ghosh, Basu and Pardo (2014).

In this context, it is useful to have a general theoretical results for the 
minimum $S$-divergence estimators beyond the discrete models, since 
there are many real life problems that can not be modeled by the discrete distributions.
The present paper fills this gap by providing a general asymptotic theory 
for the minimum $S$-divergence estimators  under continuous model families.
To avoid the complications of kernel bandwidth selection, 
we here considered the Basu-Lindsay (1994) approach of smoothed densities 
to define minimum $S^*$-divergence estimator; we have also discussed the equivalence 
of the MSDE and MSDE$^*$ under suitable assumptions 
and discussed its practical implications. All the theoretical results derived have been 
well-supported by extensive simulation study and real data examples.

\begin{acknowledgements}
The authors thank the editor and two anonymous referees for several useful suggestions and comments
that lead to an improved version of the manuscript.
\end{acknowledgements}


\begin{thebibliography}{}
\bibitem{P00}
	Pearson, K. (1900). On the criterion that a given system of deviations from the probable in
	the case of a correlated system of variables is such that it can be reasonably supposed to
	have arisen from random sampling. {\em Philosophical Magazine}, {\bf 50}, 157--175.
\end{thebibliography}


\end{document}